\documentclass[a4paper,12pt]{amsart}
\usepackage[latin1]{inputenc}
\usepackage[T1]{fontenc}
\usepackage{amsfonts,amssymb}
\usepackage{pstricks}
\usepackage{amsthm}
\usepackage[all]{xy}
\usepackage[pdftex]{graphicx}
\usepackage{graphicx}   
\newtheorem*{Theo}{Theorem}

\newtheorem{theo}{Theorem}[section]
\newtheorem{defi}{Definition}[section]
\newtheorem{prop}[theo]{Proposition}

\newtheorem{lem}[theo]{Lemma}
\newtheorem{conj}[theo]{Conjecture} 

\newtheorem{rem}[theo]{Remark}
\newtheorem{ex}{Example}

\title[Looking for a new version of Gordon's identities]{Looking for a new version of Gordon's identities}

\author{Pooneh Afsharijoo}

\begin {document}

\maketitle  
\begin{abstract}

We give a commutative algebra viewpoint on Andrews recursive formula for the partitions appearing in \textit{Gordon's identities}, which are a generalization of Rogers-Ramanujan identities. Using this approach and differential ideals we conjecture a family of partition identities which extend Gordon's identities. This family is indexed by $r\geq 2.$ We prove the conjecture for $r=2$ and $r=3.$

\end{abstract}

\footnote{{\textbf{2010 Mathematics Subject Classification.} 
05A17,12H05,13D40,13P10.}}

\maketitle
\section{INTRODUCTION}
A \textit{partition} (of length $\ell$) of a positive integer $n$ is a sequence $\Lambda:(\lambda_1\geq \cdots \geq \lambda_\ell)$ of positive integers $\lambda_i$, for $1 \leq i \leq \ell$, such that
$$\lambda_1+\cdots +\lambda_\ell=n.$$
The integers $\lambda_i$ are called the \textit{parts of the partition $\Lambda$}.
\\
The number of different partitions of $n$ is denoted by $p(n).$ By convention we set  $p(0)=1$. 
\\
A partition identity is an equality for every $n$ between the number of the partitions of an integer $n$ satisfying a certain condition $A$ and the number of those satisfying another condition $B$. This type of identity plays an important role in many areas such as number
theory, combinatorics, Lie theory, particle physics and statistical mechanics. In general it is difficult to find partition identities and to prove them. See \cite{A} for a detailed exposition of partition theory. In this article we will use commutative algebra to find and prove some partition identities.
\\
Our bridge between commutative algebra and partitions is the \textit{Hilbert-Poincar\'e series}:

Let  $ \mathbf  k$ be a field of characteristic zero. Recall that the Hilbert-Poincar\'e series of a graded $ \mathbf k-$algebra $B=\bigoplus_{i \in \mathbb{N}} B_i$ such that $\dim_{ \mathbf k}(B_i)< \infty$ is by definition the following $q-$series:

$$HP(B)=\sum_{i\in \mathbb{N}} \dim_{ \mathbf k}(B_i)q^i.$$
The generating series for the partition function $p(n)$ is given by

$$\sum_{n=0}^{\infty} p(n)q^n= \prod_{i\geq 1}\frac{1}{1-q^i}.$$

We can see that this is equal to the  Hilbert-Poincar\'e series of the graded algebra $ S=\mathbf k[x_1,x_2,\cdots]$ where the grading is given by wt.$x_i=i$ (Note that this ring is the algebra of the global sections of the space of arcs centered at the origin of the affine line (see Section $2$)).
\nolinebreak Indeed, to each monomial $x_{\alpha_m}\cdots x_{\alpha_1}$ of weight $\sum_{j=1}^{m}\alpha_j=n,$ we can associate a unique partition $(\alpha_1,\cdots,\alpha_m)$ of $n$ where $\alpha_1 \geq \cdots \geq \alpha_m.$

One important family of partitions identities is the family of \textit{Gordon's identities} (see Theorem $1$ in \cite{G}):

\begin{theo} \label{Gordon}  (\textit{Gordon's identities}). Given integers $r\geq 2$ and $1\leq i \leq r,$ let $B_{r,i}(n)$ denote the number of partitions of $n$ of the form $(b_1,\dots, b_s)$, where $b_{j}-b_{j+r-1} \geq 2$ and at most $i-1$ of the integers $b_j$ are equal to $1$. Let $A_{r,i}(n)$ denote the number of partitions of $n$ into parts $\not\equiv 0,\pm i \ ( \text{mod}.  2r+1)$. Then $A_{r,i}(n)=B_{r,i}(n)$ for all integers $n$.

\end{theo}

This is Theorem $7.5$ in [And98].
Corollary $7.9$ in [And98] gives the analytic form of this theorem, which is as follows:

\begin{theo} \label{cor} (Gordon's identities, analytic form)  For the integers $2\leq r, \ 1 \leq i \leq r,$ we have
$$\sum_{n_1,n_2,\dots n_{r-1} \geq 0} \frac{q^{N_1^2+N_2^2+\dots +N_{r-1}^2+N_i+ N_{i+1}+\dots+N_{r-1}}}{(q)_{n_1}(q)_{n_2}\dots (q)_{n_{r-1}}}=\prod_{\underset{n \not\equiv 0, \pm i (mod. 2r+1)}{n\geq1}} \frac{1}{1-q^n}.$$
Where $q$ is a variable and $N_j=n_j +n_{j+1} +\dots + n_{r-1}$ for all $1 \leq j\leq r-1 $ and $(q)_n=(1-q)(1-q^2)\cdots(1-q^n)$.
\end{theo}

The left hand side of the equality above is the generating series of $B_{r,i}(n)$ and its right hand side is the generating series of $A_{r,i}(n).$
\\
A celebrated special case of this theorem, which is known in the literature as \textit{The first Rogers-Ramanujan identity} (respectively \textit{The second Rogers-Ramanujan identity}), is when we take $r=i=2$ (respectively $r=i+1=2$). 

In this paper, we study partition identities using the relation between partitions and the graded algebras associated to an important object of algebraic geometry:  the space of arcs. We only need to consider the space of arcs of the algebraic $\mathbf{k}$-scheme defined by $(x^r)\subset \mathbf k[x]$ for any integer $r \geq 2$ (for the definition in the general case see Section $2$). This corresponds to the set $X_\infty=\{x(t)\in \mathbf k[[t]]| \ x^r(t)=0\},$ where $ \mathbf k[[t]]$ is the formal power series ring in one indeterminate $t.$ Since $x(t)\in \mathbf k[[t]]$ we can write it as $\sum_{i\in \mathbb N} x_i t^i$ and hence $x^r(t)$ is also a formal power series in $t.$ We denote the coefficients of $t^i$ in this series by $F_i;$ note that $F_i\in \mathbf k [x_0,x_1,\cdots].$ 

The space of arcs centered at the origin is obtained by setting $x_0=0$ in $X_\infty.$ Its corresponding algebra is:
$$J_\infty ^0(X)=\frac{S}{(F_{i_{|x_0=0}}|i\geq 1)}.$$
We call it the \textit{focussed arc algebra of $X$}. If we define the derivation $D$ on $S$ by  $D(x_i)=x_{i+1}$ and we denote the ideal $(x_1^r,D^1(x_1^r),D^2(x_1^r),\cdots)$ by $I_r,$ then we observe that (see Section $2$):

$$J_\infty ^0(X)  \simeq \frac{S}{I_r}.$$

The ideal $I_r$ is a differential ideal in the sense that $D(I_r) \subset I_r.$ 
%We denote the arc space of $X$ by $X_\infty$
%space correspond to the quotient of $k[x_1,x_2,\cdots]$ by the ideal 

%Given the polynomials $f_1,\cdots,f_m \in \mathbf k[x_1,\cdots,x_n],$ the associated arc space is the space corresponds to the ideal $I$ of $\mathbf k[x_{1,j},\cdots,x_{n,j}| \ j\geq 0],$ generated by the coefficients of the expansion 
%$$f_\ell(\sum_{j\geq 0}x_{1,j} t^j,\cdots,\sum_{j\geq0}x_{n,j}t^j),$$ for $1\leq \ell \leq m.$
%The ideal $I$ is graded by wt.$x_i=i.$ 
%For focussed arcs, which is when we take $x_{i,0}=0$ for $1\leq i \leq n,$ the Hilbert-Poincar\'e series of the graded algebra $\frac{\mathbf k[x_{1,j},\cdots,x_{n,j}| \ j\geq 0]}{I}$ is closely related to partitions of integers satisfying conditions depending on $I.$

%In the case where $f(x)=x^r\in \mathbf k[x],$ the ideal $I_r:=I\subset \mathbf k[x_1,x_2,\cdots]$ is a differential ideal for the derivation $D(x_i)=x_{i+1},$ in the sens that $DI_r \subset I_r.$ In fact it is generated by $x_1^r$ and all its iterated derivatives (see Section $2$).

 We use the correspondence explained above for the ring $S$ between the generating series of the partitions and the Hilbert-Poincar\'e series of the graded algebras to do the following:

%In [AB], G. E. Andrews and R. J. Baxter proved this special case. They started from the right-hand  side of the equality of the corollary above; for $r=i=2$ they denoted it by $G_2(q),$ and for $r-1=i=1$ by $G_1(q)$. Then they defined the infinite family of series $G_i(q)$ for $i\geq3$ as follow:
%\\
%$$G_i(q)=\frac{G_{i-2}(q)-G_{i-1}(q)}{q^{i-2}}$$
%\\
%They proved Rogers-Ramanujan identities, by proving the following facts which collectively are called the \textit{Empirical Hypothesis}: 
%\\For each $1 \leq i,$ 
%\begin{itemize}
%\item $G_i(q)$ contains only non-negative powers of $q,$
%\item The constant term of $G_i(q)$ is equal to $1,$
%\item $q^i$ divides $G_i(q)-1.$
%\end{itemize}

%In [LZ],  J. Lepowsky, M. Zhu denoted the product side of the equality of the above corollary by $G_i$ For each integer $1 \leq i \leq r$, and then for $j\ge1$ and $i=2, \dots, r$, they defined the formal series
%$$G_{(r-1)j+i}=\frac{G_{(r-1)(j-1)+r-i+1}-G_{(r-1)(j-1)+r-i+2}}{q^{(i-1)j}}.$$ 

%Then for all $1\leq i \leq r,$ they empirically found a recursion formula for $G_i$ and they proved Gordon's identities by generalizing an idea that have been used by G. E. Andrews and R. J. Baxter in [AB]. 

\begin{itemize}

\item[1.] 
%In \cite{LZ}, J. Lepowsky, M. Zhu gave another proof of Gordon's identities. They empirically found a recursion formula for the generating series of the partitions which are counted by $A_{r,i}(n)$ and then they generalized the idea used by G.
%E. Andrews and R. J. Baxter in \cite{AB}. 

We express the generating series of $B_{r,i}(n)$ as the Hilbert-Poincar\'e series of the quotient of $S$ by the ideal $$I_{r,i}=(x_1^i,L_{<_{\text{revlex}}}(I_r))=(x_1^i, x_j^{r-n} x_{j+1}^n| \ j \geq 1 \ \& \ \  0 \leq n \leq r-1),$$ which we determine from $J_\infty^0(X)$  and a result of C. Bruschek, H. Mourtada, J. Schepers (Proposition $5.2$ from \cite{BMS}).

Using the properties of Hilbert-Poincar\'e series  we find a recursion formula for the generating series of $B_{r,i}(n)$. Then we show that this recursion formula is equal to that found empirically by J. Lepowsky and M. Zhu to prove Gordon's identities (see \cite{LZ}). This is done in Section $3.$

\item[2.]
% In Theorem $1.6$ in \cite{AM}, we proved a theorem whose special case adds a new member to the first Rogers-Ramanujan identity. We took the algebraic variety $X$ define by $(x^2)\subset \mathbf k[x]$. computing Hilbert-Poincar\'e series of the graded algebras related to the differential ideal $I_2$ and Hilber-Poincar\'e series properties we proved this theorem.

 In the fourth section, we use the \textit{Andrews-Baxter system}  (see the proof of Theorem $7.5$ from \cite{A}) to obtain a family of identities of Rogers-Ramanujan type. A part of this family was proved in Theorem $1.6$ in \cite{AM} by the Hilbert-Poincar\'e series method. But here we use the language of partitions and we show that there is another type of partitions of $n$ whose number satisfy an extension of Andrews-Baxter system. This gives us the following theorem:

 \begin{Theo}\label{new}  (The $k$-th identity of Rogers-Ramanujan type )
 Let $n, \ m\geq0, k\geq 1$ be integers and $i=1$ or $2$. Let us denote by  $c_{2,i}^k(m,n)$ the number of partitions of $n$ of the form $(\lambda_1,\lambda_2,\cdots,\lambda_m)$, such that $\lambda_m> m+k-i$ . Let denote by $b_{2,i}^k(m,n)$  the number of partitions of $n$  of the form $(\lambda_1,\lambda_2,\cdots,\lambda_m)$ with $\lambda_m\geq k$, at most $i-1$ parts equal to $k$ and without equal or consecutive parts.
  Then $c_{2,i}^k(m,n)=b_{2,i}^k(m,n).$ 
\end{Theo}

Note that when we take $k=1$ in the theorem above we obtain a new version of the Rogers-Ramanujan identities.
%we prove a theorem which gives us a new version of the two Rogers-Ramanujan identities (see Theorem \ref {new}), but this time we prove it just by using the language of partitions.

\item[3.] In the fifth section, we conjecture a new version of Gordon's identities. To do so we use the method introduced to prove Theorem $1.6$ in \cite{AM}: We take the algebraic scheme $X$ defined by $(x^r)\subset \mathbf k[x].$ We know that (see e.g., Theorem $5.2.6$ in \cite{GP}) the Hilbert-Poincar\'e series of a homogeneous weighted ideal is equal to the Hilbert-Poincar\'e series of its leading ideal with respect to any monomial ordering. Using this theorem twice we obtain:

$$HP(\frac{S}{L_{<_{\text{revlex}}}(I_r)})=HP(\frac{S}{I_r})=HP(\frac{S}{L_{<_{\text{ lex}}}(I_r)}).$$

Note that the left hand side of this equality is the generating series of the number of  partitions of $n$ which appear on one side of Gordon's identities for $i=r$ (see \cite{BMS}).
\\
Note also that its right hand side is the generating series of the numbers of partitions of $n$ associated to the monomials of weight $n$ which appear in the graded algebra $\frac{S}{L_{<_{\text{ lex}}}(I_r)}.$
\\
 Thus, by the above equality the numbers of these partitions could give a new version of Gordon's identities for $i=r.$

In order to find this new version, we tried to find a Gr\"obner basis of the ideal $I_r$ with respect to the weighted lexicographic order. Theorem $2.2$ from \cite{AM} shows that such a Gr\"obner basis is differentially infinite in the case $r=2$ and so is complicated to compute in the general case. We could not find such a Gr\"obner basis, but  the computations give us a candidate for $L_{<_{\text{ lex}}}(I_r)$. Using this candidate we conjecture a new version of Gordon's identities for every $1\leq i \leq r$ (see Conjecture \ref{conjG} and also Conjecture \ref{conjG1} ).
The theorem above proves this conjecture for $r=2.$
 In the last four sections we give this conjecture and its analytic form (see Conjecture \ref{anal i=r}). Then we prove it also for  $r=3$.  The results presented here are part of my Ph.D. thesis \cite{Af}.

\end{itemize}
\section*{ACKNOWLEDGMENT}

I would like to express my deep gratitude to Hussein Mourtada, my Ph.D advisor, for suggesting me this project and for his constant help. I am very thankful to Bernard Teissier for his help and his corrections of the earlier version of this paper. I also would like to thank Jehanne Dousse and Fr\'ed\'eric  Jouhet for showing interest in this work and for giving me the opportunity to talk about it in their seminar. I also would like to thank the two anonymous reviewers for their useful suggestions and comments.

\section{SPACE OF ARCS}

In this section we recall the definition of the space of arcs of an algebraic scheme. %Note that an algebraic scheme defined by an ideal $I,$ is the zero set of all polynomials of the ideal $I.$ 

Let $\mathbf k$ be a field of characteristic zero and $m,n\geq 1$ be  integers. Let $X$ be an algebraic scheme defined by the ideal $(f_1,\cdots,f_m) \subset \mathbf k[x_1,\cdots,x_n].$ The \textit{arc space} of $X$, that we denote by $X_\infty$ is the set 

$$X_\infty=\{x(t)=(x_1(t)\cdots,x_n(t))\in \mathbf k[[t]]^n| \ f_\ell(x(t))=0 \ \  \text{for all}\  1\leq \ell \leq m\},$$ 

where $\mathbf k[[t]]$ denotes the formal power series ring in one variable $t$ over the field $\mathbf k.$ It has a natural structure of a $\mathbf k$-scheme. Since for each $1\leq i \leq n$ we have $x_i(t)\in \mathbf k[[t]]$, we can write it as $x_i(t)=\sum_{j=0}^\infty x_{i,j} t^j$ and we have

$$f_\ell(x(t))= F_{\ell,0} + F_{\ell,1} t +F_{\ell,2} t^2+\cdots,$$

where $F_{\ell,j} \in \mathbf k[x_{i,j} | 1 \leq i \leq n,\  0\leq j]$ for all $1\leq \ell \leq m.$ Note that $F_{\ell,0}=f_\ell(x_{1,0},x_{2,0}\cdots,x_{n,0}).$ We assume that the point $(0,\dots,0) \in X;$ hence for all $1\leq \ell \leq m$ we have $F_{\ell,0}(0,\cdots,0)=f_\ell(0,\cdots,0)=0.$

 %We are interested in the case where we have $F_{\ell,0}=0,$ which means that $(x_{1,0},x_{2,0}\cdots,x_{n,0})$ is a point of the algebraic variety $X.$ Without loss of generality we may assume that $F_{\ell,0}=0.$ 
 
 For each $i \in \{1,\cdots,n\}$ and $j \in \mathbb{N}_{>0},$ we replace $x_{i,0}$ by $0$ in $F_{\ell,j},$ and we denote the resulting polynomial by $f_{\ell,j}.$ Then we obtain a $\mathbf k$-algebra which is called \textit{the focussed arc algebra of $X$}:

$$J_\infty^0(X)=\frac{\mathbf k[x_{i,j} | 1 \leq i \leq n,\  1\leq j]}{(f_{\ell,j} | 1\leq \ell \leq m,\ 1\leq j)}.$$

The spectrum of $J_\infty^0(X)$ is the space of arcs of $X$ centered at the origin $(0,\cdots,0).$

If we give $x_{i,j}$ the weight $j,$ then $f_{\ell,j}$ is quasi-homogeneous of weight $j.$ Thus, $J_\infty^0(X)$ is a naturally graded algebra where the (usual) degree of each monomial $x_{i_1}^{j_1}\cdots x_{i_n}^{j_n}\in \mathbf k[x_1,x_2,\cdots]$ is $\sum_{k=1}^{n}j_k,$ and its weight is $\sum_{k=1}^{n}i_k j_k$.

%\begin{ex}
%Consider the algebraic variety $X$ defined by $(x^2) \subset \mathbf k[x].$ We have:

%$$X_\infty=\{ x(t) \in \mathbf k[[t]], x^2(t)=0\}.$$
%Since 

%$$x^2(t)=(x_0+x_1t+x_2t^2+\cdots)^2$$

%$$=\underbrace{x_0^2}_{F_0}+\underbrace{2x_0x_1}_{F_1}t+\underbrace{(2x_0x_2+x_1^2)}_{F_2}t^2+\underbrace{(2x_0x_3+2x_1x_2)}_{F_3}t^3+\cdots,$$

%we have  $f_0=f_1=0, f_2=x_1^2, f_3=2x_1x_2, \cdots$ and so

%$$J_\infty^0(X)=\frac{\mathbf k[x_1,x_2,\cdots]}{(f_1,f_2,\cdots)}.$$

%$x^r(t)=(x_0+x_1t+x_2t^2+x_3t^3+\cdots)^r=(\frac{y_0}{0!}+\frac{y_1}{1!}t+\frac{y_2}{2!}t^2+\frac{y_3}{3!}t^3+\cdots)^r$
%$=y_0^2+2y_0y_1t+(y_0y_2+y_1)^2t^2+(\frac{y_0}{3}+y_1y_2)t^3+\cdots$

%$=\frac{y_0^2}{0!}+\frac{2y_0y_1}{1!}t+\frac{2y_0y_2+2y_1^2}{2!}t^2+\frac{2y_0y_3+6y_1y_2}{3!}t^3+\cdots$

%$= \frac{y_0^2}{0!}+\frac{D^1(y_0^2)}{1!}+\frac{D^2(y_0^2)}{2!}+\frac{D^3(y_0^2)}{3!}+\cdots$

%\end{ex}

In the \textit{weighted lexicographic order}, we say that $x^\alpha=x_1^ {{\alpha}_1}\cdots x_n^{{\alpha}_n}$ is less than $x^\beta=x_1^ {{\beta}_1}\cdots x_n^{{\beta}_n}$ if and only if wt.$x^\alpha < \text {wt.} x^\beta$ or  wt.$x^\alpha =\text{wt.}x^\beta$ and there exists $j\geq 1$ such that $\alpha_1=\beta_1, \cdots, \alpha_{j-1}=\beta_{j-1}, \alpha_j<\beta_j$.

The \textit{weighted reverse lexicographic order} consists also in comparing first the weights,  in case of equality of the weights, we have $x^\alpha < x^\beta$ if and only if there exists $j\geq 1$ such that $\alpha_l=\beta_l$ for all $l>j$ and $\alpha_j>\beta_j$.

We fix a monomial ordering $>$ on $\mathbf k[x_1,x_2,\cdots].$ For  $f\in \mathbf k[x_1,x_2,\cdots]$ the \textit{leading monomial} of $f$ is its largest monomial with respect to $>.$ The \textit{leading ideal} of an ideal $I \subset \mathbf k[x_1,x_2,\cdots]$ is the ideal generated by the leading monomials of the polynomials in $I.$ Note that in general it is not equal to the ideal generated by the leading monomials of generators of $I.$ A set of nonzero polynomials of $I$ whose leading monomials with respect to $>$ generate the leading ideal of $I$ is called a \textit{Gr\"obner basis} of $I$ with respect to $>.$

The following example of space of arcs is fundamental for our work:

\begin{ex} \label{arc} Consider the algebraic scheme $X$ defined by $(f) \subset \mathbf k[y],$ where $f=y^r$ for some integer $r \geq 2$ (for reasons that will appear below, we change the name of the variable). So the arc space of this algebraic scheme is the set 

$$X_\infty=\{ y(t) \in \mathbf k[[t]], y^r(t)=0, y(t)=\sum_{j=0}^{\infty}y_j t^j\}.$$ 

We have $y^r(t)=\Big(\sum_{j=0}^{\infty} y_jt^j\Big)^r=\sum_{j=0}^{\infty} F_jt^j.$
\\
\\
Let $x_j=y_j j!$ and let $D$ be the derivation on $\mathbf k[x_0,x_1,\cdots]$ defined by $D(x_i)=x_{i+1}.$ We have

%Making change of variable $y_j$ by $\frac{x_j}{j!}$  and defining the derivation $D$ on $\mathbf k[x_0,x_1,\cdots]$ as $D(x_i):=x_{i+1}$, we can give a differential structure to the focussed arc algebra of $X$ (see in \cite{M}). We have:
$y^r(t)=(y_0+y_1t+y_2t^2+y_3t^3+\cdots)^r=(\frac{x_0}{0!}+\frac{x_1}{1!}t+\frac{x_2}{2!}t^2+\frac{x_3}{3!}t^3+\cdots)^r$
\\
\\
$= \frac{x_0^r}{0!}+\frac{D^1(x_0^r)}{1!}t+\frac{D^2(x_0^r)}{2!}t^2+\frac{D^3(x_0^r)}{3!}t^3+\cdots$
\\
\\
Let $f_i=F_i\mid_{x_0=0}.$ So we have (see Proposition 2.1 in \cite{M1}):
 
$$\frac{\mathbf k[y_0,y_1,\cdots]}{(F_0,F_1,\cdots)} \simeq \frac{\mathbf k[x_0,x_1,\cdots]}{(x_0^r,D^1(x_0^r),D^2(x_0^r),\cdots)},$$

where the $i$-th derivation $D^i$ is recursively defined by $D^1(g)=D(g)$ and $D^i(g)=D(D^{i-1}(g))$ for all $g \in \mathbf k[x_0,x_1,\cdots].$ We also obtain:

$$\frac{\mathbf k[y_1,y_2,\cdots]}{(f_0,f_1,\cdots)} \simeq \frac{\mathbf k[x_1,x_2,\cdots]}{(x_1^r,D^1(x_1^r),D^2(x_1^r),\cdots)}.$$
%If we define the derivation $D$ on $\mathbf k[x_1,x_2,\cdots]$ as $D(x_i):=x_{i+1}$, we have $F_0=x_0^r$ and $F_\ell=D(F_{\ell-1})$ for all $1\leq \ell$. 

Thus, the differential ideal  $(x_1^r,D^1(x_1^r),D^2(x_1^r),\cdots)$ is the ideal defining the focussed arc algebra of $X$ in $S=\mathbf k[x_1,x_2,\dots]$.
We denote this differential ideal by $I_r=[x_1^r].$
\\
  In \cite{BMS}, C. Bruschek, H. Mourtada, J. Schepers  proved that the leading ideal of $I_r$ with respect to the reverse lexicographical order is as follows:
$$L(I_r)=I_{r,r}=(x_j^{r-n} x_{j+1}^n| \ j \geq 1 \ \& \ \  0 \leq n \leq r-1).$$

From now on, we consider the focussed arc algebra after this change of variable.
\end{ex}

For $1\leq i \leq r$ define the ideal $I_{r,i}=(x_1^i,I_{r,r}).$ Let us consider the graded algebra $\frac{S}{I_{r,i}}$. We have:
 $$dim_{\mathbf k}\Big( \frac{S}{I_{r,i}}\Big)_j =\dim_{\mathbf k} \frac{S_j}{(I_{r,i})_j} \leq  \dim_{\mathbf k} S_j=p(j)< \infty.$$ 
So the Hilbert-Poincar\'e series of  $\frac{S}{I_{r,i}}$ exists and has the following form:

$$HP_{\frac{S}{I_{r,i}}}(q)=\sum_{j\in \mathbf{N}}\dim_{\mathbf k} \Big( \frac{S}{I_{r,i}}\Big)_j q^j.$$

To each monomial $x_{i_m}\cdots x_{i_1}\in S$ we can associate a partition $(i_1,\cdots,i_m)$ where $i_1\geq \cdots \geq i_m.$ Since $I_{r,i}$ is generated by $x_1^i, x_1^{i-1}x_2, \cdots , x_1 x_2^{r-1}$ and the monomials of the form $x_j^{r-n}x_{j+1}^n$ where $j\geq 2$ and $0\leq n \leq r-1,$ computing the Hilbert-Poincar\'e series of the graded algebra $\frac{S}{I_{r,i}}$ is equivalent to counting the partitions which are counted by $B_{r,i}(n)$ in Gordon's identities (see Theorem \ref{Gordon}). i.e.,
 $$HP_{\frac{S}{I_{r,i}}}(q)=\sum_{n\geq 0} B_{r,i}(n)q^n.$$   
\\
Note that an important property of the Hilbert-Poincar\'e series is that, if $E \subset S$ is a homogeneous ideal and $f\in S$ is a homogeneous polynomial of degree $d$ then we have the following exact sequence (see Lemma 5.2.2 in [GP])

$$0\longrightarrow \frac{S}{(E:f)}(-d) \longrightarrow \frac{S}{E} \longrightarrow \frac{S}{(E,f)} \longrightarrow 0,$$
%to avvali g ro mifrestim be fg
where $(E:f)=\{g \in S | \ fg \in E\}.$ So we have (see Corollary 6.2 in \cite{BMS})
\begin{equation} \label{eq:1}
HP_{\frac{S}{E}}(q)=q^d HP_{\frac{S}{(E:f)}}(q)+HP_{\frac{S}{(E,f)}}(q).
\end{equation}

We will use this equation many times in this paper. 

\section{A RECURSION FORMULA FOR $HP_{\frac{S}{I_{r,i}}}(q)$ VIA      COMMUTATIVE ALGEBRA}

In this section we prove a recursion formula producing formal power series which converge to $HP_{\frac{S}{I_{r,i}}}(q)$ in the $q$-adic topology. This proves Gordon's identities. To do this we need some notations.
\\
\\
For each integer $k\geq 1$ denote $\mathbf k[x_k,x_{k+1},\cdots]$ by $S_k.$ We shall use the following ideals of $S_k:$
\\
$$J_k=(x_i^{r-n} x_{i+1}^n, \ i \geq k,\  0 \leq n \leq r-1),$$

$$J_k^l=(x_k^l,x_k^{l-1}x_{k+1}^{r-l+1},x_k^{l-2}x_{k+1}^{r-l+2},\cdots, x_k x_{k+1}^{r-1}, J_{k+1}),$$
\\
where  $1\leq l \leq r.$ In this section we will denote  the Hilbert-Poincar\'e series $HP_{\frac{S_k}{J_k}}(q)$ by $H^k$ and the Hilbert-Poincar\'e series $HP_{\frac{S_k}{J_k^l}}(q)$ by $H_l^k$. Note that $H^k_1=H^{k+1}$ and $H^k_r=H^k.$ Note also that:
$$HP_{\frac{S}{I_{r,i}}}(q)=HP_{\frac{S_1}{J_1^i}}(q)=H_i^1.$$ 
\\

To construct the recursion formula for $H_i^1$, we use the following two lemmas:

\begin{lem} \label{ser1} With the notations introduced above, for $1\leq l \leq r$ and $k\geq 1$, we have
$$H_l^k=\sum_{j=1}^{l} q^{(l-j)k}H_{r-l+j}^{k+1}.$$

%$=q^{(l-1)k}H_{r-l+1}^{k+1}+q^{(l-2)k}H_{r-l+2}^{k+1}+q^{(l-3)k}H_{r-l+3}^{k+1}+\cdots+q^kH_{r-1}^{k+1}+H^{k+1}.$$
\end{lem}
\begin{proof}  For simplicity we omit $q$ and we use $HP(\frac{A}{I})$ Instead of  $HP_{\frac{A}{I}}(q).$ Using Equation (\ref {eq:1}) we have:

$$H_l^k=q^k HP(\frac{S_k}{(J_k^l:x_k)})+HP(\frac{S_k}{(J_k^l,x_k)}) $$

%$$=q^k HP(\frac{S_k}{((x_k^l,x_k^{l-1}x_{k+1}^{r-l+1},x_k^{l-2}x_{k+1}^{r-l+2},\cdots, x_k x_{k+1}^{r-1}, J_{k+1}):x_k)})$$

%$$+HP(\frac{S_k}{(x_k,x_k^l,x_k^{l-1}x_{k+1}^{r-l+1},x_k^{l-2}x_{k+1}^{r-l+2},\cdots, x_k x_{k+1}^{r-1}, J_{k+1})}).$$

$$=q^k  HP(\frac{S_k}{(x_k^{l-1},x_k^{l-2} x_{k+1}^{r-l+1},x_k^{l-3}x_{k+1}^{r-l+2},\cdots, x_k x_{k+1}^{r-2}, J_{k+1}^{r-1})})+H^{k+1}.$$

We continue in this way and we use repetitively Equation (\ref{eq:1}), so we obtain

$$H_l^k=q^k \Big(q^k HP(\frac{S_k}{\big((x_k^{l-1},x_k^{l-2} x_{k+1}^{r-l+1},x_k^{l-3}x_{k+1}^{r-l+2},\cdots, x_k x_{k+1}^{r-2}, J_{k+1}^{r-1}):x_k\big)})$$

$$+HP(\frac{S_k}{(x_k,x_k^{l-1},x_k^{l-2} x_{k+1}^{r-l+1},x_k^{l-3}x_{k+1}^{r-l+2},\cdots, x_k x_{k+1}^{r-2}, J_{k+1}^{r-1})})\Big)+H^{k+1}$$

$$=q^{2k}  HP(\frac{S_k}{(x_k^{l-2},x_k^{l-3} x_{k+1}^{r-l+1},x_k^{l-4}x_{k+1}^{r-l+2},\cdots,  x_k x_{k+1}^{r-3}, J_{k+1}^{r-2})})+q^k H_{r-1}^{k+1} +H^{k+1}$$

$$=\cdots =\sum_{j=1}^{l} q^{(l-j)k}H_{r-l+j}^{k+1}.$$

\end{proof}

Using the previous lemma we can give a formula for $H^k$:

%\begin{lem} \label{ser2}For $k\geq1$, we have
%$$H^k=q^{(r-1)k} H^{k+2}+q^{(r-2)k}H_2^{k+1} +q^{(r-3)k}H_3^{k+1}+\cdots+q^kH_{r-1}^{k+1}+H^{k+1}.$$

%\end{lem}

%\begin{proof}
 %The proof is similar to the proof of Lemma \ref{ser1}, using repetitively Equation (\ref{eq:1}) this time starting with $H^k$ .
 
 % For simplicity we omit $q$ and we use $HP(\frac{A}{I})$ Instead of  $HP_{\frac{A}{I}}(q).$ Using Equation (\ref{eq:1}) repetitively,  we have
%$$H^k=q^kHP(\frac{S_k}{(J_k:x_k )})+HP(\frac{S_k}{(x_k, J_k)})$$

%$$=q^k HP(\frac{S_k}{(x_k^{r-1},x_k^{r-2}x_{k+1},x_k^{r-3} x_{k+1}^2,\cdots,x_k x_{k+1}^{r-2},J_{k+1}^{r-1})})+H^{k+1}$$

%$$=q^k \Big( q^k HP(\frac{S_k}{\big((x_k^{r-1},x_k^{r-2}x_{k+1},x_k^{r-3} x_{k+1}^2,\cdots,x_k x_{k+1}^{r-2},J_{k+1}^{r-1}): x_k \big)}$$

%$$ + HP(\frac{S_k}{(x_k,x_k^{r-1},x_k^{r-2}x_{k+1},x_k^{r-3} x_{k+1}^2,\cdots,x_k x_{k+1}^{r-2},J_{k+1}^{r-1})})\Big) +H^{k+1}$$

%$$=q^{2k} HP(\frac{S_k}{(x_k^{r-2},x_k^{r-3}x_{k+1},x_k^{r-4}x_{k+1}^2,\cdots ,x_k x_{k+1}^{r-3},J_{k+1}^{r-2})})$$

%$$+q^k H_{r-1}^{k+1}+H^{k+1}=\cdots$$

%$$=q^{(r-1)k} H^{k+2}+q^{(r-2)k}H_2^{k+1} +q^{(r-3)k}H_3^{k+1}+\cdots+q^kH_{r-1}^{k+1}+H^{k+1}.$$
%\end{proof}
We are now ready to give the recursive formula for $H_i^1:$

\begin{prop} \label{H1} For all integers $r \geq 2,\ 1\leq i \leq r$ we have the following recursion formula:

%$$H_i^1= B_{i,1,(r-1)(d-1)+1} H^d + B_{i,2,(r-1)(d-1)+2} H_{r-1}^d$$
%$$+\cdots +B_{i,r-1,(r-1)(d-1)+r-1}H_2^d+B_{i,r,(r-1)(d-1)+r}H^{d+1},$$

$$H_i^1=\sum_{j=1}^{r} B_{i,j,(r-1)(d-1)+j}H_{r-j+1}^d.$$
\\
where $d\geq 3$ and the $B_{i,j,k} \in \mathbf k[[q]]$ satisfy the following recursion formula for $1\leq l\leq r$
%$$B_{i,j,(r-1)(d-1)+j}= q^{(j-1)(d-1)} (B_{i,1,(r-1)(d-2)+1}+B_{i,2,(r-1)(d-2)+2}$$
%$$+\cdots +B_{i,(r-j+1),(r-1)(d-2)+r-j+1}).$$

$$B_{i,j,(r-1)(d-1)+j}=q^{(j-1)(d-1)}\sum_{k=1}^{r-j+1}B_{i,k,(r-1)(d-2)+k}.$$
With the following initial conditions:

$B_{i,j,2r+j-2}=\begin{cases}
q^{2(j-1)} (1+q+\cdots+q^{i-1}) &\text{ if\ }\  1\leq j \leq r-i+1 \\
q^{2(j-1)} (1+q+\cdots+q^{r-j}) &\text{ if\ }\  r-i+2 \leq j \leq r.
\end{cases}$
\end{prop}

\begin{proof}
The proof is by induction on $d$. Assume $d=3$, by Lemma \ref{ser1} for $k=1$ and $l=i,$ we have
 %$$H_i^1=q^{(i-1)}H_{r-i+1}^2+q^{(i-2)}H_{r-i+2}^2+q^{(i-3)}H_{r-i+3}^2+\cdots+qH_{r-1}^2+H^2.$$
 $$H_i^1=\sum_{j=1}^{i} q^{i-j}H^2_{r-i+j}.$$
 Now, using Lemma \ref{ser1}we replace $H_l^2,$ for $r-i+1 \leq l \leq r,$ in the  equation above:
 %  $H_i^1=q^{(i-1)}(q^{2(r-i)}H_i^3+q^{2(r-i-1)}H_{i+1}^3+q^{2(r-i-2)}H_{i+2}^3+\cdots+q^2 H_{r-1}^3+H^3)$
%   $+q^{(i-2)}(q^{2(r-i+1)}H_{i-1}^3+q^{2(r-i)}H_i^3+q^{2(r-i-1)}H_{i+1}^3+\cdots+q^2 H_{r-1}^3+H^3)$
 %  $+q^{(i-3)}(q^{2(r-i+2)}H_{i-2}^3+q^{2(r-i+1)}H_{i-1}^3+q^{2(r-i)}H_i^3+\cdots+q^2 H_{r-1}^3+H^3)$
 %  $\vdots$
%   $+q(q^{2(r-2)}H_2^3+q^{2(r-3)}H_3^3+q^{2(r-4)}H_4^3+\cdots+q^2 H_{r-1}^3+H^3)$
  % $+(q^{2(r-1)}H^4+q^{2(r-2)}H_2^3+q^{2(r-3)}H_3^3+\cdots+q^2 H_{r-1}^3+H^3).$
   $$H_i^1=\sum_{j=1}^iq^{i-j}\sum_{k=1}^{r-i+j}q^{2(r-i+j-k)}H_{i-j+k}^3.$$
Factoring out  $H_l^3$ for $1\leq l \leq r$ proves our formula for $d=3.$

%\begin{equation*}
	%\begin{split}
		%H_i^1&=(1+q+\cdots+q^{(i-1)})H^3\\
		%&+q^2(1+q+\cdots+q^{(i-1)})H_{r-1}^3\\
		%&+q^4(1+q+\cdots+q^{(i-1)})H_{r-2}^3\\
		%&\vdots\\
		%&+q^{2(r-i)}(1+q+\cdots+q^{(i-1)})H_i^3\\
		%&+q^{2(r-i+1)}(1+q+\cdots+q^{(i-2)})H_{i-1}^3\\
		%&+q^{2(r-i+2)}(1+q+\cdots+q^{(i-3)})H_{i-2}^3\\
		%&\vdots\\
		%&+q^{2(r-2)}(1+q)H_2^3\\
		%&q^{2(r-1)}H^4.
	%\end{split}
%\end{equation*}
Let us now assume that the formula is true for $d\leq m$ and prove it for $d=m+1$. By the induction hypothesis for $d=m$ we obtain this expression for $H_i^1:$
% $$B_{i,1,(r-1)(m-1)+1} H^m + B_{i,2,(r-1)(m-1)+2} H_{r-1}^m+\cdots +B_{i,r-1,(r-1)(m-1)+r-1}H_2^m+B_{i,r,(r-1)(m-1)+r}H^{m+1}.$$
 $$H_i^1=\sum_{j=1}^{r} B_{i,j,(r-1)(m-1)+j}H_{r-j+1}^m.$$
 
By Lemma \ref{ser1} we have:
%$H_i^1= B_{i,1,(r-1)(m-1)+1} \Big(q^{(r-1)m} H^{m+2}+q^{(r-2)m}H_2^{m+1} +\cdots+q^m H_{r-1}^{m+1}+H^{m+1}\Big)$
%$+B_{i,2,(r-1)(m-1)+2} 
%\Big(q^{(r-2)m}H_{2}^{m+1}+q^{(r-3)m}H_{3}^{m+1}+\cdots+q^m H_{r-1}^{m+1}+H^{m+1}\Big)$
%$\vdots$
%$+B_{i,r-1,(r-1)(m-1)+r-1} \Big(q^{m}H_{r-1}^{m+1}+H^{m+1} \Big)$
%$+B_{i,r,(r-1)(m-1)+r}H^{m+1}.$
$$H_i^1=\sum_{j=1}^rB_{i,j,(r-1)(m-1)+j}\sum_{k=1}^{r-j+1}q^{m(r-j+1-k)}H_{j-1+k}^{m+1}.$$
We rewrite now the  equation above in another way by factoring out  $H_l^{m+1}$ for $1\leq l \leq r:$
%$H^1_i=\Big (B_{i,1,(r-1)(m-1)+1}+B_{i,2,(r-1)(m-1)+2}+\cdots +B_{i,r,(r-1)(m-1)+r}\Big) H^{m+1}$
%$+q^m  \Big(B_{i,1,(r-1)(m-1)+1}+B_{i,2,(r-1)(m-1)+2}+\cdots +B_{i,r-1,(r-1)(m-1)+r-1}\Big) H_{r-1}^{m+1}$
%$+q^{2m}  \Big(B_{i,1,(r-1)(m-1)+1}+B_{i,2,(r-1)(m-1)+2}+\cdots +B_{i,r-2,(r-1)(m-1)+r-2}\Big) H_{r-2}^{m+1}$
%$\vdots$
%$+q^{(r-2)m}  \Big(B_{i,1,(r-1)(m-1)+1}+B_{i,2,(r-1)(m-1)+2}\Big) H_2^{m+1}$
%$+q^{(r-1)m}  B_{i,1,(r-1)(m-1)+1} H^{m+2}.$
$$H_{i}^1=\sum_{l=1}^{r}q^{m(r-l)}\sum_{j=1}^l B_{i,j,(r-1)(m-1)+j}H_l^{m+1},$$
 which by our notations and the recursion formula of $B_{i,j,k},$ this is the same as:
%$H_i^1= B_{i,1,(r-1)m+1} H^{m+1} + B_{i,2,(r-1)m+2} H_{r-1}^{m+1}$
%$+\cdots+B_{i,r-1,(r-1)m+r-1}H_2^{m+1}+B_{i,r,(r-1)m+r}H^{m+2}.$
$$H_i^1=\sum_{l=1}^rB_{i,r-l+1,m(r-1)+r-l+1}H_l^{m+1}=\sum_{j=1}^{r} B_{i,j,m(r-1)+j}H_{r-j+1}^{m+1}.$$
\end{proof}

Now fix an integer $r \ge 2$. For each $l=1,\dots , r$ define
$$G_l =\prod_{\underset{n \not\equiv 0,\pm (r+1-l) (\text{mod}. 2r+1)}{n\geq1}} \frac{1}{1-q^n}.$$
Note that $G_l$ is the product side of the equation in the analytic form of Gordon's identities (see Theorem \ref{cor} ) where $i=r+1-l$.
We want to show that for $i=r+1-l$, $G_l $ is equal to $H_i^1$, which proves Gordon's identities. To do this, we show that $G_l$ and $H_i^1$ are both limits for the $q$-adic topology of the same sequence of polynomials in $q$.
\\
We now use a recursion formula of J. Lepowsky, M. Zhu in \cite{LZ}.
\\
\\
For $j\ge1$ and $i=2, \dots, r$, define recursively the formal power series
$$G_{(r-1)j+i}=\frac{G_{(r-1)(j-1)+r-i+1}-G_{(r-1)(j-1)+r-i+2}}{q^{(i-1)j}}.$$ 

So we have 
\begin{equation}\label{LZ}
 G_{(r-1)j-i+2}=q^{(i-1)j} G_{(r-1)j+i} +G_{(r-1)j-i+3}.
\end{equation}

\begin{prop} \label{G1} (\textsc{J. Lepowsky, M. Zhu in \cite{LZ}}) For all integer $1\leq l \leq r$, we have the following recursion formula:
 
$$G_l= \sum_{j=1}^{r}A_{l,j,(r-1)d+j}G_{(r-1)d+j},$$
where $ d\geq2$ and $A_{l,j,k} \in \mathbf k[[q]]$ satisfy the following recursion formula for $1\leq j\leq r$
$$A_{l,j,(r-1)d+j}=q^{(j-1)d}\sum_{k=1}^{r-j+1}A_{l,k,(r-1)(d-1)+k}.$$
With the following initial condition:
\\
\\
$A_{l,j,2r+j-2}=\begin{cases}
q^{2(j-1)} (1+q+\cdots+q^{r-l}) &\text{ if\ }\  1\leq j \leq l \\
q^{2(j-1)} (1+q+\cdots+q^{r-j}) &\text{ if\ }\  l+1 \leq j \leq r.
\end{cases}$
\end{prop}

%So far, we have constructed recursion formulas for $H_i^1$ and $G_l.$ 
Note that in \cite{LZ} the polynomial $A_{l,j,(r-1)d+j}$ is denoted by $_l h^{(d)}_j$. We will now show that if $i=r-l+1,$ the coefficients of the two formulas of Propositions \ref{H1} and \ref{G1} are equal:
\begin{prop} \label{coef} With the notations used in this section, for all $d\geq 2$ and $1\leq m \leq r$, we have 
$$A_{l,m,(r-1)d+m}=B_{i,m,(r-1)d+m},$$ 
where $ 2\leq r,\ 1\leq i \leq r$ and $l=r-i+1$.

\end{prop}

\begin{proof}The proof is by induction on $d$. Note that by Proposition \ref{G1}, we have
\\ 
$A_{l,m,2(r-1)+m}=A_{l,m,2r+m-2}=\begin{cases}
q^{2(m-1)} (1+q+\cdots+q^{r-l}) &\text{ if\ }\  1\leq m \leq l \\
q^{2(m-1)} (1+q+\cdots+q^{r-m}) &\text{ if\ }\  l+1 \leq m \leq r.
\end{cases}$
\\
Replacing $l$ by $r-i+1$ we obtain
\\
$A_{l,m,2(r-1)+m}=\begin{cases}
q^{2(m-1)} (1+q+\cdots+q^{i-1}) &\text{ if\ }\  1\leq m \leq r-i+1 \\
q^{2(m-1)} (1+q+\cdots+q^{r-m}) &\text{ if\ }\  r-i+2 \leq m \leq r.
\end{cases}$
\\
This is equal to $B_{i,m,2(r-1)+m}$ by Proposition \ref {H1}.
Now assume that the equation is true for $d-1$.
Again by Propositions \ref{H1}, for all $1\leq m \leq r$,  we have
%$$B_{l,m,(r-1)d+m}= q^{(m-1)d} (B_{l,1,(r-1)(d-1)+1}+B_{l,2,(r-1)(d-1)+2}
%+\cdots +B_{l,(r-m+1),(r-1)(d-1)+r-l+1}).$$
$$B_{l,m,(r-1)d+m}=q^{(m-1)d}\sum_{k=1}^{r-m+1}B_{l,k,(r-1)(d-1)+k}.$$
By the induction hypothesis we obtain
%$$B_{l,m,(r-1)d+m}= q^{(m-1)d} (A_{i,1,(r-1)(d-1)+1}+A_{i,2,(r-1)(d-1)+2}+\cdots +A_{i,(r-m+1),(r-1)(d-1)+r-m+1}).$$
$$B_{l,m,(r-1)d+m}=q^{(m-1)d}\sum_{k=1}^{r-m+1}A_{l,k,(r-1)(d-1)+k}.$$

By Proposition \ref{G1}, the right hand side of the  equation above is equal to $A_{i,m,(r-1)d+m}$.

\end{proof}
Now we are ready to prove the main theorem of this section, which gives another proof for \textit{Gordon's identities}.
\begin{theo} \label{equ} With the notations used in this section, we have
$$G_l=H_i^1,$$
where $2\leq r,\ 1\leq i \leq r$ and $l=r-i+1$.

\end{theo}

 \begin{proof} We denote the limit of a sequence of formal power series $a_i \in \mathbf k[[q]]$ in the $q$-adic topology (if it exists) by $\lim a_i.$ By Proposition \ref{G1}, since the power of $q$ in $A_{l,m,(r-1)d+m}$ is greater than or equal to $(m-1)d,$
 %$\text{ord} A_{l,m,(r-1)d+m} \geq q^{(m-1)d}$ for all $2 \leq m \leq r$, 
 it is immediate that $ \lim_{d\to+\infty} A_{l,m,(r-1)d+l}$ exists for all $1 \leq m \leq r$; in fact for all $2\leq m \leq r$ we have
 
 $$ \lim_{d\to+\infty} A_{l,m,(r-1)d+m}= \lim_{d\to+\infty} q^{(m-1)d} \sum_{k=1}^{r-m+1}A_{l,k,(r-1)(d-1)+k}=0,$$
 
 and so $$ G_l= \lim_{d\to+\infty} A_{l,1,(r-1)d+1} G_{(r-1)d+1}.$$
 %$$ \lim_{d\to+\infty} q^{(m-1)d} (A_{l,1,(r-1)(d-1)+1}+\cdots +A_{l,(r-m+1),(r-1)(d-1)+r-m+1})=0.$$

  Theorem $2.1$ of \cite{LZ} implies that $ G_{(r-1)d+i}$ is a formal power series with the constant term equal to $1$ and that $ G_{(r-1)d+i}-1$ is divisible by $q^{d+1}$ if $1\leq i \leq r-1$ and by $q^{d+2}$ if $i=r$ (this is the Empirical Hypothesis of Lepowsky and Zhu). Thus  $ \lim_{d\to+\infty} G_{(r-1)d+1}=1 $ and we have:
 $$ G_l= \lim_{d\to+\infty} A_{l,1,(r-1)d+1} G_{(r-1)d+1}=\lim_{d\to+\infty} A_{l,1,(r-1)d+1}  .$$

 Let us denote $ \lim_{d\to+\infty} A_{l,1,(r-1)d+1} $ by $A_{l,1,\infty}$. So we have 
 $G_l=A_{l,1,\infty}  .$
 
 On the other hand, in the same way as above, by Proposition \ref{H1}  for all $2 \leq m \leq r$ we have:
$$\lim_{d\to+\infty} B_{i,m,(r-1)(d-1)+m}=\lim_{d\to+\infty} q^{(m-1)(d-1)} \sum_{k=1}^{r-m+1}B_{i,k,(r-1)(d-1)+k}.$$
  %$$\lim_{d\to+\infty} q^{(m-1)(d-1)} (B_{i,1,(r-1)(d-2)+1}+\cdots +B_{i,(r-m+1),(r-1)(d-2)+r-m+1}(q))=0.$$
  So we have 

 $$ H_i^1 =\lim_{d\to+\infty} B_{i,1,(r-1)(d-1)+1} H^d.$$
   Note that $\lim_{d\to+\infty} H^d=1.$ Because the zero-th homogeneous component is isomorphic to $\mathbf{k}$ and hence is of dimension $1.$ The homogeneous component of degree $i,$ for $1\leq i <d,$ is zero since there is no monomials of degree between $1$ and $d.$ So:

   % by the Hilbert-Poincar\'e series properties $H^d$ is equal to:
   
  $$H^d=1+q^d \alpha(q),$$
  where $\alpha(q)\in \mathbf{k}[|q|].$ Thus, $\lim_{d\to+\infty} H^d=1$ and we have:
  
  $$ H_i^1=\lim_{d\to+\infty} B_{i,1,(r-1)(d-1)+1}.$$
 
 If we denote  $ \lim_{d\to+\infty} B_{i,1,(r-1)(d-1)+1}$ by $B_{i,1,\infty},$ we have
 $$  H_i^1 =B_{i,1,\infty}.$$
 
 By Proposition $3.5$ for all $d\geq2$ and $1\leq m \leq r$, we have $A_{l,m,(r-1)d+m}=B_{i,m,(r-1)d+m}$ where $l=r-i+1$. Hence $A_{l,1,\infty}=B_{i,1,\infty}$. So we have
 $$G_l=A_{l,1,\infty}=B_{i,1,\infty}=H_i^1.$$
 \end{proof}

 \section{THE $K$-TH ROGERS-RAMANUJAN TYPE  IDENTITY }
 
If we take $r=i=2$ (respectively $r=i+1=2$) in Gordon's Identities we obtain a special case which is called \textit{the first Rogers-Ramanujan identity} (respectively \textit{the second Rogers-Ramanujan identity}). In \cite{AM} by using the Hilbert-Poincar\'e series properties we proved a theorem whose special case gives a new version of the first Rogers-Ramanujan identity.

In this section we prove the following theorem, which gives us not only this version of the first Rogers-Ramanujan identity, but also a new version of the second Rogers-Ramanujan identity, but this time by using the language of partitions:

 \begin{theo}\label{new} 
 Let $n, \ m\geq0$ and $k\geq 1$ be integers and let $c_{2,i}^k(m,n)$ denote the number of partitions of $n$ of the form $(\lambda_1,\lambda_2,\cdots,\lambda_m)$, such that $\lambda_m> m+k-i$ for $i=1,2.$ Let $b_{2,i}^k(m,n)$ denote the number of partitions of $n$  of the form $(\lambda_1,\lambda_2,\cdots,\lambda_m)$ with $\lambda_m\geq k$, at most $i-1$ parts equal to $k$ and without equal or consecutive parts.
  Then $c_{2,i}(m,n)=b_{2,i}(m,n).$ 
\end{theo}

\begin{proof}
 We prove that the $c_{2,i}^k(m,n)$ satisfy an extention of Andrew's system (see the proof of Theorem  $7.5$ from \cite{A}). In other words we want to prove:
\\
$c_{2,i}^k(m,n)=\begin{cases}{\label{1}}
1 &\text{ if } m=n=0 \\
0 &\text{ if } m\leq 0 \text{ or } n\leq 0 \text{ but } (m,n) \neq (0,0);
\end{cases}$

$c_{2,2}^k(m,n)-c_{2,1}^k(m,n)=c_{2,1}^k(m-1,n-m-k+1);$
 
$c_{2,1}^k(m,n) =c_{2,2}^k(m,n-m).$
\\
\\
Note that $0$ has only one partition whose length is zero (the empty set). A negative number has no partition, and a positive number has no partition of non positive length. So the first equation is true.
\\
\\
For the second one, note that the left hand side of this equation counts the number of partitions of $n$ of the form  $(\lambda_1,\lambda_2,\cdots,\lambda_{m-1},m+k-1)$. If we delete $m+k-1$ from this partition, we obtain a partition of $n-m-k+1$ with exactly $m-1$ parts $(\lambda_1,\lambda_2,\cdots,\lambda_{m-1})$ such that the last part $,\lambda_{m-1},$ is at least equal to $m+k-1$. This defines a one-to-one correspondence between the partitions counted by $c_{2,2}^k(m,n)-c_{2,1}^k(m,n)$ and those counted by $c_{2,1}^k(m-1,n-m-k+1).$
\\
\\
For the last equation, we will transform each partition of $n$ of the form $(\lambda_1,\lambda_2,\cdots,\lambda_m)$ with exactly $m$ parts such that $\lambda_m >m+k-1$ by subtracting 1 from each part. Since $\lambda_m \geq m+k$, obviously $\lambda_m -1 \geq m+k-1$. So by this transformation we obtain the partitions of $n-m$ with exactly $m$ parts such that the smallest part is at least equal to $m+k-1.$ Thus, to each partition counted by $c_{2,1}^k(m,n)$ we associated a unique partition which is counted by $c_{2,2}^k(m,n-m).$ Obviously this transformation also is a bijection between the partitions counted by each side of this equation, which proves the last equation. 

So far we proved that  $c_{2,i}^k(m,n)$ satisfy the above system. Note that using the same method as Andrews, one can show that $b_{2,i}^k(m,n)$  are \textit {uniquely} determined by this system of equations (see the proof of Theorem 7.5 in \cite{A}). Therefore $c_{2,i}^k (m,n)=b_{2,i}^k (m,n)$ for all $m$ and $n$ with $ i =1, 2.$

%Now let $b_{2,i}(m,n)$ denote the number of partitions $(b_1,b_2,\cdots,b_s)$ of $n$ with exactly $m$ parts such that $b_j-b_{j+1}\geq 2$ and for the integer $1 \leq i \leq 2$ at most $i-1$ of the $b_j$ equal $1.$ 

\end{proof}

\begin{rem}
Let $C_{2,i}^k(n)=\sum_{m \geq 0} c_{2,i}^k (m,n),$ for $i=1,2.$ 
Note that by Theorem \ref{new} we have:

$$C_{2,i}^k(n)=\sum_{m \geq 0} c_{2,i}^k (m,n)=\sum_{m \geq 0} b_{2,i}^k (m,n).$$
For $k=1$ the right hand side of the equality above is equal to $B_{2,i}(n).$ So in this case we obtain a new version of Rogers-Ramanujan identities.
%\linebreak where $B_{2,i}(n)$ and $A_{2,i}(n)$ are the same as in Gordons identities.
\\
We proved also that even if we fix a length for the partitions, the number of partitions of $n$ counted by $C^1_{2,i}(n)$ will always be equal to the number of those partitions counted by $B_{2,i}(n).$ This is not true for $A_{2,i}(n)$ and $B_{2,i}(n)$ and so it is not true in general for Gordon's identities.
\end{rem}

\section{A NEW CLASSIFICATION OF PARTS OF A PARTITION AND GORDON'S IDENTITIES}

In Example \ref{arc} we took the algebraic scheme $X$ defined by $(x^r)\subset \mathbf k[x]$ for some integer $r\geq 2,$ and we saw that the ideal which defines its focussed arc algebra is the differential ideal $I_r=[x_1^r].$ As we mentioned in the introduction, in order to find a new version of Gordon's identities we wanted to find the Gr\"obner basis of the ideal $I_r$ with respect to the weighted lexicographic order. Since such a Gr\"obner basis is differentially infinite in the case $r=2$ (see Theorem $2.2$ from \cite{AM}) and so it is complicated to compute in the general case, we could not find it in general. But the computations suggested a candidate for  $L_{<_{\text{ lex}}}(I_r)$, let us denote this candidate by $I'_{r,r}.$

 This candidate gives us a conjecture for a new version of Gordon's identities (see Conjecture \ref{conjG} and also Conjecture \ref{conjG1}). This conjecture claims that the Hilbert-Poincar\'e series associated to our candidate monomial ideal is the generating series for the number of partitions satisfying certain conditions related to \textit{the new parts of a partition} defined below:

\begin{defi} \label{def} Given an integer $r\geq 2,$ for $1 \leq i \leq r$ we define $(i,\ell)$-new part of $\Lambda:(\lambda_1,\cdots,\lambda_m),$ as follows:

$$p_{i,\ell}(\Lambda)=\begin{cases}
\lambda_m &\text{ if } \ell=1 \\
\lambda_{m-\sum_{j=1}^{\ell-1}p_{i,j}(\Lambda)} &\text{ if }  2\leq \ell \leq i \\
\lambda_{m+\ell-i-\sum_{j=1}^{\ell-1}p_{i,j}(\Lambda)} &\text{ if } i <\ell \leq r-1;
\end{cases}$$ 
where $\lambda_j=0$ for $j\leq 0,$ and if $p_{i,\ell}(\Lambda)=0$ then $p_{i,j}(\Lambda)=0$ for $j> \ell.$ We denote the number of all non zero  $(i,\ell)$-new part of $\Lambda$ by $N_i(\Lambda).$
\end{defi}

In the other words we define $p_{i,\ell}(\Lambda)$ recursively as follows:
\begin{itemize}
\item $p_{i,1}(\Lambda)=$The smallest part of $\Lambda$;

\item For $2\leq \ell \leq i,$ we define the $(i,\ell)$-new part of $\Lambda$ as its $(\sum_{j=1}^{\ell-1}p_{i,j}(\Lambda)+1)$-th part counting from the right;
% agar vojood nadasht, barabar ba sefr dar nazar migirim.
\item For $i+1\leq \ell \leq r$  we define the $(i,\ell)$-new part of $\Lambda$ as its $(\sum_{j=1}^{\ell-1}p_{i,j}(\Lambda)+i-\ell+1)$-th part counting from the right.
\end{itemize}

%and if $p_{i,\ell}(\Lambda)=0$ then $p_{i,j}(\Lambda)=0$ for all $j> \ell.$
% dar vaghe in shart ham dg lazem nis, chon agar ye ja p_{i,\ell}(\Lambda)=0, yani joze p
Note that $1\leq N_i(\Lambda) \leq r-1.$ Let us look at an easy example to become more familiar with these new parts.
 \begin{ex}
 Take the partition $\Lambda:(4,4,3,2,2,2)$ of $17.$ For $r=i=4$ we have
 \\
 \\
$p_{4,1}(\Lambda)=2;$
\\
$p_{4,2}(\Lambda)=$The third part of $\Lambda$ counting from the right $=2;$
\\
$p_{4,3}(\Lambda)=$ The fifth part of $\Lambda$ counting from the right $=4;$
\\
So $N_4(\Lambda)=3.$
\end{ex}

%Let denote by $C_{r,i}(n),$ the number of partitions of an integer $n$ which have at most $i-1$ parts equal to $1,$ and whose length and $(i,\ell)$-new parts satisfy one of the following conditions:

%\begin{itemize}
%\item The partitions with $N_i(\Lambda)<r-1$;
%\item The partitions with $N_i(\Lambda)=r-1$ whose length is at most equal to the sum of its $(i,\ell)$-new parts plus $i-r$.
%\end{itemize}

We are now ready to state our conjecture:

%\begin{ex} 
% Take the partition $\Lambda:(4,4,3,2,2,2)$ of $17.$ For $r=i=4$ we have
 
%$p_{4,1}(\Lambda)=\lambda_m=\lambda_6=2;$

%$p_{4,2}(\Lambda)=\lambda_{m-p_{4,1}(\Lambda)}=\lambda_{6-2}=2;$

%$p_{4,3}(\Lambda)=\lambda_{m-p_{4,1}(\Lambda)-p_{4,2}(\Lambda)}=\lambda_{6-2-2}=4.$

%\end{ex}

\begin{conj}\label{conjG} (A new version of Gordon's identities) For an integer $n,$ let $C_{r,i}(n)$ denote the number of partitions of $n$ of the form $\Lambda:(\lambda_1,\cdots,\lambda_s),$ such that at most $i-1$ of the parts $\lambda_j$ are equal to $1$ and either $N_i(\Lambda)<r-1,$ or $N_i(\Lambda)=r-1$ and  $s\leq \sum_{j=1}^{r-1} p_{i,j}(\Lambda)-(r-i)$. Then $C_{r,i}(n)=B_{r,i}(n)=A_{r,i}(n),$ where $B_{r,i}(n),A_{r,i}(n)$ are the same as in Gordon's identities.
\end{conj}

\begin{rem}\label{conj1} If we define $p_{i,r}(\Lambda):=\lambda_{m+r-i-\sum_{j=1}^{r-1}p_{i,j}(\lambda)}$ then Conjecture \ref{conjG} can be expressed as follows:
\\
\\
 For an integer $n,$ let $C_{r,i}(n)$ denote the number of partitions of $n$ whose $(i,\ell)$-new part is equal to zero for some $1 \leq \ell \leq r $. Then $C_{r,i}(n)=B_{r,i}(n)=A_{r,i}(n),$ where $B_{r,i}(n),A_{r,i}(n)$ are the same as in Gordon's identities.

\end{rem}
Note that Theorem \ref {new} proves this conjecture for $r=2.$

In order to prove this conjecture we defined the ideal $I'_{r,i}=(x_1^i,I'_{r,r})$ and we proved  that the Hilbert-Poincar\'e series of the graded algebra $\frac{S}{I'_{r,i}}$ is equal to the generating series of $C_{r,i}(n)$ in Conjecture \ref{conjG} (see Proposition \ref{I'}). The problem is that we could not prove the equality between this $q$-series and the generating series of $B_{r,i}(n)$ (or $A_{r,i}(n)$). 
%that this Hilbert-Poincar\'e series is equal to one of the formulas which appear in the analytic form of Gordon's identities.

% we first tried to show that for all $1\leq i \leq r,$ the generating series of the numbers of partitions counted by $C_{r,i}(n)$ is equal to the generating series of the partitions counted by $B_{r,i}(n).$ i.e.,

%$$\sum_{n\geq 0} C_{r,i}(n) q^n = \sum_{n\geq 0} B_{r,i}(n) q^n.$$

%To do so, we took a field $\mathbf k$ of characteristic zero and we considered a weighted polynomial ring $S=\mathbf k[x_1,x_2,\cdots]$ where $\text{wt}. x_i=i.$ Then we defined an ideal $I'_{r,i}$ such that the Hilbert-Poincar\'e series of the graded algebra $\frac{S}{I'_{r,i}}$ is the generating series of the number of partitions of an integers $n\geq 1$ counted by $C_{r,i}(n)$ in Conjecture \ref{conjG}. Then we tried to prove that this Hilbert-Poincar\'e series is equal to one of the formulas which appear in the analytic form of the Gordon's identities.
Let us now introduce $I'_{r,i}.$
 To do so we define $r$ blocks of increasing positive integers with the following property: 
 \\
 The first block contains only one integer. For $2\leq j \leq i,$ the number of integers which appear in the $j-$th block is equal to the last number of the previous block. For $i+1\leq j \leq r,$  the number of integers which appear in the $j-$th block is equal to the last number of the previous block minus one.  i.e.,
 
 $$\underbrace{n_{1,1}}_{\text {The first block}}\leq \underbrace{ n_{2,1} \leq \cdots \leq n_{2,n_{1,1}} }_{\text {The second block}}\leq \underbrace{ n_{3,1} \leq \cdots \leq n_{3,n_{2,n_{1,1}}}}_{\text{The third block}}\leq \cdots.$$
 
In order to simplify notations, for $1\leq j \leq r$ we introduce:

$$f(j)= \begin{cases}
1 &\text{ if } j=1 \\
n_{j-1,f(j-1)} &\text{ if } 2\leq j \leq i \\
n_{j-1,f(j-1)}-1 &\text{ if } i+1 \leq j \leq r+1.
\end{cases}$$

So we are considering the following $r$ blocks of positive integers:

  $$\underbrace{n_{1,1}}_{\text {The first block}}\leq \underbrace{ n_{2,1} \leq \cdots \leq n_{2,f(2)} }_{\text {The second block}}\leq \underbrace{ n_{3,1} \leq \cdots \leq n_{3,f(3)}}_{\text{The third block}}\leq \cdots \leq \underbrace{ n_{r,1} \leq \cdots \leq n_{r,f(r)} }_{\text {The $r$-th block}}.$$

Now let  $I'_{r,i}$ be the ideal generated by  monomials of the form:

$$x_{n_{1,1}}x_{n_{2,1}}\cdots x_{n_{2,f(2)}} x_{n_{3,1}} \cdots x_{n_{3,f(3)}} \cdots  x_{n_{r,1}}\cdots x_{n_{r,f(r)}}.$$

Recall that to each monomial $x_{\alpha_1}\cdots x_{\alpha_k} \in S$  we can associate a partition $(\alpha_k,\cdots, \alpha_1)$ where $\alpha_1\leq \cdots \leq \alpha_k.$
 Thus to each generator
$$x_{n_{1,1}}x_{n_{2,1}}\cdots x_{n_{2,f(2)}} x_{n_{3,1}} \cdots x_{n_{3,f(3)}} \cdots  x_{n_{r,1}}\cdots x_{n_{r,f(r)}}$$
of $I'_{r,i}$, we can associate the following partition:  
$$(n_{r,f(r)},\cdots,n_{r,1},\cdots,n_{2,f(2)},\cdots,n_{2,1},n_{1,1}).$$
If as usual, we denote this partition by 
$\Lambda:(\lambda_1,\cdots,\lambda_m)$
 then $m=\sum_{j=1}^{r} f(j)$ and $\lambda_{m-\sum_{j=1}^{s} f(j)+1}=n_{s,f(s)}$ for all $1\leq s \leq r.$

\begin{prop} \label{I'} For integers $r\geq 2$ and $1\leq i \leq r$ we have:
$$HP_{\frac{S}{I'_{r,i}}}(q)=\sum_{n\geq 0} C_{r,i}(n)q^n.$$
\end{prop}

\begin{proof} For each partition $\Lambda:(\lambda_1,\cdots,\lambda_m)$ we take:

\begin{itemize}
\item  $n_{1,1}:=\lambda_m;$
\item $n_{s,k}:=\lambda_{m-\sum_{j=1}^{s-1}f(j)-k+1}$ where $2\leq s$ and $1\leq k \leq f(s)-1;$
\item  $n_{s,f(s)}:=\lambda_{m-\sum_{j=1}^{s}f(j)+1}$  where $2\leq s;$

\end{itemize}
where $\lambda_j=0$ for $j\leq0.$
Note that by induction on $j$ one can show that:
$$p_{i,j}(\Lambda)= \begin{cases}
f(j+1) &\text{ if } 1\leq j \leq i-1 \\
f(j+1)+1 &\text{ if } i \leq j \leq r.
\end{cases}$$

 Consider now a partition $\Lambda:(\lambda_1,\cdots,\lambda_m)$ which is counted by $C_{r,i}(n)$. By Remark \ref{conj1} this means that $p_{i,\ell}(\Lambda)=0$ for some $1\leq \ell \leq r.$ Therefore, by definition of the new parts, we have $p_{i,j}(\Lambda)=0$ for all $j\geq \ell.$ Thus:

 \begin{equation*}
\begin{aligned}
\Lambda \ \text{is counted by} \ C_{r,i}(n) &
\iff \exists 1\leq \ell \leq r, \ \forall \ell \leq j \leq r, \ p_{i,j}(\Lambda)=0;\\ 
      &\iff  \exists 1\leq \ell \leq r, \ \forall \ell\leq j \leq r, \ \begin{cases}
f(j+1) &\text{ if } 1\leq j \leq i-1 \\
f(j+1)+1 &\text{ if } i \leq j \leq r.
\end{cases}=0; \\ 
&\iff  \exists 1\leq \ell \leq r, \ \forall  \ell\leq j \leq r, \ \begin{cases}
n_{j,f(j)} &\text{ if } 1\leq j \leq i-1 \\
(n_{j,f(j)}-1)+1 &\text{ if } i \leq j \leq r.
\end{cases}=0;\\
&\iff  \exists 1\leq \ell \leq r, \ \forall  \ell\leq j \leq r, \ n_{j,f(j)}=0;\\
&\iff  \exists 1\leq \ell \leq r, \ x_{\Lambda}=x_{n_{1,1}}x_{n_{2,1}}\cdots x_{n_{2,f(2)}}\cdots x_{n_{\ell,1}}\cdots x_{n_{\ell,s}}\\
&\ \ \ \ \ \ \ \ \ \ \ \ \ \ \ \ \ \ \ \ \ \ \ \ \ \ \ \ \ \text{where} \ 1\leq s \leq f(\ell)-1;\\
&\iff  x_{\Lambda} \notin I'_{r,i};\\
& \iff x_{\Lambda} \in \frac{S}{I'_{r,i}} .\\
\end{aligned}
\end{equation*}

%So in this case $\mu$ is a partition which is counted by $C_{r,i}(\mu_1+\cdots+\mu_m)$.
 This proves the equality between the Hilbert-Poincar\'e series of the graded algebra $\frac{S}{I'_{r,i}}$ and the generating series of $C_{r,i}(n)$.
\end{proof}

%So the Hilbert-Poincar\'e series of $\frac{S}{I'_{r,i}}$ is the generating series of $C_{r,i}(n)$ and so the following Conjecture is equivalent to  Conjecture \ref {conjG}:
%\begin{conj}\label{Conj 2} For any integer $r\geq3,$ we have
%$$HP_{\frac{S}{I'_{r,i}}}(q)=\sum_{n_1,n_2,\dots n_{r-1} \geq 0} \frac{q^{N_1^2+N_2^2+\dots +N_{r-1}^2+N_i+ N_{i+1}+\dots+N_{r-1}}}{(q)_{n_1}(q)_{n_2}\dots (q)_{n_{r-1}}}=\prod_{\underset{n \not\equiv 0, \pm i (mod. 2r+1)}{n\geq1}} \frac{1}{1-q^n}.$$

%\end{conj}
%We tried to prove this conjecture by computing the Hilbert-Poincar\'e series of the graded algebra $\frac{S}{I'_{r,i}}$ and using the Hilbert-Poincar\'e series properties. But unfortunately we have not succeeded. 

In the last section of this paper we compute $HP_{\frac{S}{I'_{r,r}}}$ and using the previous proposition, we state the analytic form of Conjecture \ref{conjG} for the case $i=r$.

%For the integers $3\leq r$ and $1\leq i \leq r,$ let $I_{r,i}$ be the following ideal:
%$$(x_1^i,x_1^{i-1}x_2^{r-i+1},\cdots,x_1x_2^{r-1},x_j^{r-\ell} x_{j+1}^{\ell}| \ 2\leq j, \ 0\leq \ell %\leq r-1).$$
%Note that $I_{r,i}=(x_1^i,I_{r,r})$ and $I'_{r,i}=(x_1^i,I'_{r,r}).$

\section{A NEW VERSION OF GORDON'S IDENTITIES FOR THE CASE $r=3$}

In this section we prove Conjecture \ref {conjG} for $r=3$ by using the language of partitions:

\begin{theo}\label{r3}  Given integers $n$ and $1 \leq i \leq 3,$ let $C_{3,i}(n)$ denote the number of partitions of $n$ whose $(i,\ell)$-new part is equal to zero for some $1\leq \ell \leq 3$. Then $C_{3,i}(n)=B_{3,i}(n)=A_{3,i}(n),$ where $B_{3,i}(n),A_{3,i}(n)$ are as in Gordon's identities.
\end{theo}

\begin{proof}
  Let $c_{3,i}(m,n)$ denote the number of partitions which are counted by $C_{3,i}(n)$ and with exactly $m$ parts. We are going to prove that $c_{3,i}(m,n)$ satisfy Andrew's system of equations. This means that we have:
  
$c_{3,i}(m,n)=\begin{cases}
1 &\text{ if } m=n=0 \\
0 &\text{ if } m\leq 0 \text{ or } n\leq 0 \text{ but } (m,n) \neq (0,0);
\end{cases}$

$c_{3,3}(m,n)-c_{3,2}(m,n)=c_{3,1}(m-2,n-m);$
 
$c_{3,2}(m,n)-c_{3,1}(m,n)=c_{3,2}(m-1,n-m);$
 
$c_{3,1}(m,n) =c_{3,3         }(m,n-m).$
\\
\\
For the first equation see the proof of Theorem \ref{new}. In order to prove the other equations, for each one we define a bijective transformation between the partitions which are counted by each side. 

To prove the second one, note that the left hand side of this equation counts the number of partitions of $n$ of the form  $\Lambda:(\lambda_1,\lambda_2,\cdots,\lambda_m)$ such that $m>\lambda_m$ and $\lambda_m+\lambda_{m-\lambda_m}=m$. We transform such a partition by removing $\lambda_m$ and $\lambda_{m-\lambda_m}$ from this partition. We obtain a partition of $n-m$ with exactly $m-2$ parts: 
%Note also that the only partition with these two properties and with two parts equal to $1$ is $\Lambda:(1,1).$ We send it to the empty set. 
%If $\Lambda:(\lambda_1,\lambda_2,\cdots,\lambda_m)\neq (1,1),$ we remove $\lambda_m$ and $\lambda_{m-\lambda_m}$ from this partition and  we obtain a partition of $n-m$ with exactly $m-2$ parts: 
\begin{itemize}

\item $\mu:(\mu_1,\cdots,\mu_{m-2})=(\lambda_1,\cdots,\lambda_{m-2}),$ if $\lambda_m=1;$
\item $\mu:(\mu_1,\cdots,\mu_{m-2})=(\lambda_1,\cdots,\lambda_{m-\lambda_m-1},\lambda_{m-\lambda_m+1},\cdots,\lambda_{m-1}),$ if $\lambda_m>1.$
\end{itemize}

If $\lambda_m=1$ then on the one hand  $\mu_{m-2}=\lambda_{m-2}\geq 2$ and on the other hand $\mu_{m-2}=\lambda_{m-2} \geq \lambda_{m-1}=m-1$. 

If $\lambda_m >1$ then $\mu_{m-2}=\lambda_{m-1}\geq \lambda_m>1.$ So $\mu_{m-2}>1$. On the one hand $\mu_{m-2}=\lambda_{m-1}\leq \lambda _{m-\lambda_{m}}=m-\lambda_m \leq m-2,.$ On the other hand we have 

$\mu_{m-2}+\mu_{m-1-\mu_{m-2}}=\lambda_{m-1}+\mu_{m-1-\lambda_{m-1}}\geq$
$\lambda_{m-1}+\mu_{m-1-\lambda_{m}} = \lambda_{m-1}+\lambda_{m-1-\lambda_{m}}\geq$

$\lambda_{m}+\lambda_{m-\lambda_{m}}=m$

So $\mu$ is a partition which is by definition counted by $c_{3,1}(m-2,n-m).$ 

%Now we want to prove that this transformation of partitions from $\Lambda$ to $\mu$ is injective. Assume that the partitions $\Lambda:(\lambda_1,\lambda_2,\cdots,\lambda_m)$ and $\Lambda':(\lambda'_1,\lambda'_2,\cdots,\lambda'_m)$  transform to $\mu:(\mu_1,\cdots,\mu_{m-2}).$ If $\mu_{m-2}\geq m-1$ we have $\lambda_i=\lambda'_i=\mu_i$ for all $1\leq i \leq m-2.$ If $\mu_{m-2}\leq m-2,$ assume that $\lambda_m< \lambda'_m,$ since $\lambda_m+\lambda_{m-\lambda_m}=m=\lambda'_m+\lambda'_{m-\lambda'_m},$ hence $\lambda_{m-\lambda_m} > \lambda'_{m-\lambda'_m}.$ So we have 

%$\mu_{m-1-\lambda_m}=\lambda_{m-1-\lambda_m}\geq \lambda_{m-\lambda_m} > \lambda'_{m-\lambda'_m} \geq \lambda'_{m+1-\lambda'_m}=\mu_{m-\lambda'_m}.$

%Since $\mu$ is a partition $m-1- \lambda_m < m-\lambda'_m$ and so $\lambda_m \geq \lambda'_m.$ But this is a contradiction (we assumed that  $\lambda_m < \lambda'_m$). So the transformation from $\Lambda$ to $\mu$ is an injection.

Let us now prove that the transformation from $\Lambda$ to $\mu$ is bijective. To do so let $\mu:(\mu_1,\cdots,\mu_{m-2})$ be a partition which is counted by $c_{3,1}(m-2,n-m).$ By definition $\mu_{m-2}>1$ and either $m-2<\mu_{m-2},$ or $1\leq \mu_{m-2} \leq m-2$ and  $\mu_{m-2}+\mu_{m-1-\mu_{m-2}}\geq m.$ 

If $\mu_{m-2}> m-2,$ we take $\Lambda:(\lambda_1,\cdots,\lambda_m)=(\mu_1.\cdots,\mu_{m-2},m-1,1),$ which is a partition of $n$ with $m$ parts and $\lambda_m+\lambda_{m-\lambda_m}=1+(m-1)=m.$ 

If $ \mu_{m-2}\leq m-2$, we will show that there exists a unique positive integer $2\leq k\leq \mu_{m-2}$ such that $\Lambda:(\mu_1,\cdots,\mu_{m-k-1},m-k,\mu_{m-k},\cdots,\mu_{m-2},k)$ is a partition of $n.$ Then we have $\lambda_m+\lambda_{m-\lambda_m}=k+(m-k)=m.$

In order to find such an integer $k,$ take the set 
$$A=\{1 \leq a\leq \mu_{m-2} | m-a > \mu_{m-a-1}\}.$$ 
%Note that if $ a \leq \mu_{m-2}$ we have $m-a\geq m-\mu_{m-2}\geq 2.$ So we can consider $\mu_{m-a-1}.$
We have $1 \in A$ and so $A$ is not empty. Let $k'=\max{A},$ since $\mu_{m-2}+\mu_{m-1-\mu_{m-2}}\geq m,$ 
%the integer $k'$ cannot be equal to $\mu_{m-2}$ and so 
we have $k'+1\leq \mu_{m-2}.$ Since $k'=\max{A},$   
 on the one hand $m-k'>\mu_{m-k'-1},$ and on the other hand $k'+1 \notin A$ and we have $m-k'-1\leq \mu_{m-k'-2}.$ Thus $k=k'+1$ is the integer that we look for. 

We now prove that $k$ is uniquely determined. Suppose that there exist two positive integers $k_1< k_2$ such that $\mu$ transforms to $\Lambda_i: (\mu_1,\cdots,\mu_{m-k_i-1},m-k_i,\mu_{m-k_i},\cdots,\mu_{m-2},k_i)$ for $i=1,2.$ Then we have $m-k_1-1 \geq m-k_2$ and:

$\mu_{m-k_1-1}\leq \mu_{m-k_2} \leq m-k_2 < m-k_1,$ which is a contradiction. So we can transform $\Lambda$ to $\mu$ which means that this transformation is also surjective and the third equality holds.

% So we can transform $\Lambda$ to $\mu$ which means that this transformation is also surjective and the third equality holds.

Let us now prove that $c_{3,2}(m,n)-c_{3,1}(m,n)=c_{3,2}(m-1,n-m).$ Note that the left hand side of this equation counts the number of partitions of $n$ of the  form $\Lambda:(\lambda_1,\lambda_2,\cdots,\lambda_m)$ such that:
\begin{itemize}
\item $\lambda_m=1, \lambda_{m-1}\geq m,$ 
or
\item $1<\lambda_m<m$ and $  \lambda_m+\lambda_{m-\lambda_m+1} \leq m+1 \leq \lambda_m+\lambda_{m-\lambda_m}.$
\end{itemize} 
We divide the set of all such partitions in to three disjoint subsets. 
Then we send these subsets by three different bijections to three disjoint sets whose union is the set of all partitions which are counted by $c_{3,2}(m-1,n-m).$ These partitions are the partitions of $n-m$ of the form $\mu:(\mu_1, \cdots, \mu_{m-1})$ with at most one part equal to $1$ such that either $\mu_{m-1}>m-2,$ or $1 \leq \mu_{m-1} \leq m-2$ and $\mu_{m-1}+\mu_{m-1-\mu_{m-1}}\geq m.$ To do so we will define the following three bijections:
\begin{itemize}

\item[bijection 1.] Sends the set of partitions $\Lambda:(\lambda_1,\cdots,\lambda_m)$ with $\begin{cases}
\lambda_m=1 \\
\lambda_{m-1}\geq m ;
\end{cases}$ to the set of partitions $\mu:(\mu_1,\cdots,\mu_{m-1})$ with
$\begin{cases}
\mu_{m-1}> m-2 ;
\end{cases}$
\item[bijection 2.] Sends the set of partitions $\Lambda:(\lambda_1,\cdots,\lambda_m)$ with 

$\begin{cases}
\ 1 < \lambda_m <m \\
\lambda_m+\lambda_{m-\lambda_m+1} = m+1 \leq \lambda_m+\lambda_{m-\lambda_m}; 
\end{cases}$ to the set of partitions

 $\mu:(\mu_1,\cdots,\mu_{m-1})$ with
$\begin{cases}
1\leq \mu_{m-1} \leq m-2\\
\mu_{m-1}+\mu_{m-{\mu_{m-1}}}\leq m-1 < \mu_{m-1}+\mu_{m-1-\mu_{m-1}};
\end{cases}$
\item[bijection 3.] Sends the set of partitions $\Lambda:(\lambda_1,\cdots,\lambda_m)$ with 

$\begin{cases}
\ 1 < \lambda_m <m \\
\lambda_m+\lambda_{m-\lambda_m+1} < m+1 \leq \lambda_m+\lambda_{m-\lambda_m}; 
\end{cases}$

to the set of partitions $\mu:(\mu_1,\cdots,\mu_{m-1})$ with
$\begin{cases}
2\leq \mu_{m-1} \leq m-2\\
\mu_{m-1}+\mu_{m-{\mu_{m-1}}}\geq m.
\end{cases}$
\end{itemize}
In order to define the first bijection, we transform $\Lambda$ to $\mu:(\lambda_1-1,\cdots,\lambda_{m-1}-1).$ So $\mu$ is a partition of $n-m$ of length $m-1$ and we have $\mu_{m-1}=\lambda_{m-1}-1\geq m-1.$ 

%Note that if $\mu_{m-1}=\mu_{m-2}=1$ we have $\lambda_{m-1}=\lambda_{m-2}=2.$ On the other hand we have $1=\mu_{m-1}\geq m-1$ so $m\leq 2$ which is a contradiction, because we supposed that  $\Lambda$ has at least $3$ parts $\lambda_m, \lambda_{m-1}$ and $\lambda_{m-2}.$ So $\mu$ cannot have more than one part equal to $1.$
Clearly the transformation from $\Lambda$ to $\mu$ is a bijection.

%defines a one-to-one correspondence between the partitions of $n$ of the form  $\Lambda:(\lambda_1,\cdots,\lambda_m)$ such that $\lambda_m=1, \lambda_{m-1}\geq m,$ and the partitions of $n-m$ of the form $\mu:(\mu_1,\cdots,\mu_{m-1})$ such that $\mu_{m-1}\geq m-1.$ 

In order to define the second bijection, we send $\Lambda$ to: 
$$\mu:(\lambda_1,\cdots,\lambda_{m-\lambda_m},\lambda_{m+2-\lambda_m}-1,\cdots,\lambda_m-1),$$
which is a partition of $n-m$ of length $m-1$. Since $\lambda_m \neq 1, \lambda_m+\lambda_{m-\lambda_m+1} = m+1$ and $\lambda_{m+2-\lambda_m}\leq \lambda_{m-\lambda_m},$ on the one hand we have
\\
$1\leq \mu_{m-1}=\lambda_m-1=m-\lambda_{m+1-\lambda_m}\leq m-\lambda_m \leq m-2 .$ 
\\
So $1\leq \mu_{m-1}\leq m-2.$ On the other hand we have  
\begin{itemize}
\item $\mu_{m-1}+\mu_{m-\mu_{m-1}}=\lambda_m-1+\mu_{m+1-\lambda_m}=\lambda_m+\lambda_{m+2-\lambda_m}-2$
\\
$\leq \lambda_m+\lambda_{m+1-\lambda_m}-2=m-1;$
\item $\mu_{m-1}+\mu_{m-1-\mu_{m-1}}=\lambda_m-1+\mu_{m-\lambda_m}=\lambda_m-1+\lambda_{m-\lambda_m}\geq m.$
\end{itemize}
%Note that if $\mu_{m-1}=\mu_{m-2}=1,$ then on the one hand, by definition of the transformation from $\Lambda$ to $\mu$, the partition $\Lambda$ has at least $3$ parts, $\lambda_m,\lambda_{m-1}$ and $\lambda_{m-2}$, which means that $m \geq 3,$ and on the other hand, since $\mu_{m-1}+\mu_{m-1-\mu_{m-1}}\geq m,$ we have $1+1\geq m$ which is a contradiction. So $\mu$ can have at most one part equal to $1$. 
%So we have transformed $\Lambda$ into a partition $\mu$ of $n-m$ which is counted by $c_{3,2}(m-1,n-m).$ 

Let us now prove that this transformation is injective. Suppose that $\Lambda:(\lambda_1,\cdots,\lambda_m)$ and $\Lambda':(\lambda'_1,\cdots,\lambda'_m)$ both transform to $\mu:(\mu_1,\cdots,\mu_{m-2})$ such that $\lambda_m < \lambda'_m$. Since $\lambda_m+\lambda_{m+1-\lambda_m}=\lambda'_m+\lambda'_{m+1-\lambda'_m}=m+1,$ hence $\lambda_{m+1-\lambda_m} > \lambda'_{m+1-\lambda'_m}.$ So we have:
$$\mu_{m-\lambda_m}=\lambda_{m-\lambda_m} \geq \lambda_{m+1-\lambda_m} >  \lambda'_{m+1-\lambda'_m}\geq  \lambda'_{m+2-\lambda'_m}=\mu_{m+1-\lambda'_m}+1>\mu_{m+1-\lambda'_m}.$$
By definition of a partition we have $m-\lambda_m<m+1-\lambda'_m,$ and so $\lambda'_m\leq \lambda_m,$ which is a contradiction.

We now prove that the transformation from $\Lambda$ to $\mu$ is surjective. Let $\mu:(\mu_1,\cdots,\mu_{m-1})$ be a partition of $n-m$ such that:

$$\begin{cases}
1\leq \mu_{m-1} \leq m-2\\
\mu_{m-1}+\mu_{m-{\mu_{m-1}}}\leq m-1 < \mu_{m-1}+\mu_{m-1-\mu_{m-1}}.
\end{cases}$$
We take $\Lambda:(\mu_1,\cdots,\mu_{m-1-\mu_{m-1}},m-\mu_{m-1},\mu_{m-\mu_{m-1}}+1,\cdots,\mu_{m-1}+1)$ which is a partition of $n$ of length $m$ and we have:
 \begin{itemize}
 \item $\lambda_m=\mu_{m-1}+1\geq 2;$
 \item $\lambda_m+\lambda_{m-\lambda_m+1}=\mu_{m-1}+1+\lambda_{m-\mu_{m-1}}=\mu_{m-1}+1+m-\mu_{m-1}=m+1;$
\item $\lambda_m+\lambda_{m-\lambda_m}=\mu_{m-1}+1+\lambda_{m-1-\mu_{m-1}}=\mu_{m-1}+\mu_{m-1-\mu_{m-1}}+1\geq m+1.$
 \end{itemize}
By definition of the transformation from $\Lambda$ to $\mu,$ clearly $\Lambda$ goes to $\mu.$ This finishes the proof of the existence of the second bijection. 
\\
For the last one, let us take a partition $\Lambda=(\lambda_1,\cdots,\lambda_m)$ of $n$ such that 
$$\begin{cases}
1 \neq \lambda_m<m \\
\lambda_m+\lambda_{m-\lambda_m+1} < m+1 \leq \lambda_m+\lambda_{m-\lambda_m}.
\end{cases}$$
Send such a $\Lambda$ to $\mu:(\lambda_1-1,\cdots,\lambda_{m-\lambda_m}-1,\lambda_{m+1-\lambda_m},\cdots,\lambda_{m-1})$ which is a partition of $n-m$ of length $m-1$ and we have:
\begin{itemize}
\item $2\leq \lambda_m \leq \lambda_{m-1}=\mu_{m-1};$
\item $\mu_{m-1}=\lambda_{m-1}\leq \lambda_{m+1-\lambda_m}<m-\lambda_m+1 \leq m-1;$
\item $\mu_{m-1}+\mu_{m-{\mu_{m-1}}}=\lambda_{m-1}+\mu_{m-\lambda_{m-1}}\geq \lambda_m+\mu_{m-\lambda_m}=\lambda_m+\lambda_{m-\lambda_m}-1 \geq m.$
\end{itemize}
%We now prove that this transformation is injective. Suppose that $\Lambda:(\lambda_1,\cdots,\lambda_m)$ and $\Lambda':(\lambda'_1,\cdots,\lambda'_m)$ both transform to $\mu:(\mu_1,\cdots,\mu_{m-1})$ and suppose that $\lambda_m<\lambda'_m.$ so $m-\lambda_m \geq m-\lambda'_m+1$ and we have

%$m+1 \leq \lambda_m+\lambda_{m-\lambda_m}=\lambda_m+\mu_{m-\lambda_m}+1 \leq \lambda_m+\mu_{m-\lambda'_m+1}+1=\lambda_m+\lambda'_{m-\lambda'_m+1}+1<\lambda'_m+\lambda'_{m-\lambda'_m+1}+1 \leq m+1,$

%which is a contradiction. So the transformation from $\Lambda$ to $\mu$ is injective.
To prove the bijectivity of this transformation, let $\mu:(\mu_1,\cdots,\mu_{m-1})$ be a partition of $n-m$ such that:
$$\begin{cases}
2\leq \mu_{m-1} \leq m-2\\
\mu_{m-1}+\mu_{m-{\mu_{m-1}}}\geq m.
\end{cases}$$
We prove that there exists a unique positive integer $2\leq k \leq \mu_{m-1}$ such that $\Lambda:(\mu_1+1,\cdots,\mu_{m-k}+1,\mu_{m-k+1},\cdots,\mu_{m-1},k)$ is a partition of $n$ of length $m$ and $k+\mu_{m-k+1}<m+1\leq k+1+\mu_{m-k}.$

If $m\leq \mu_{m-2}+2,$ we take $k=2.$ We have
$2+\mu_{m-1}\leq m< m+1 \leq \mu_{m-2}+3.$

If $m> \mu_{m-2}+2,$ define the set 
$$B=\{2\leq b\leq\mu_{m-1}| m>\mu_{m-b}+b\}.$$
% Note that if $b\in \mu_{m-1}$ then we have $m-b\geq m- \mu_{m-1}\geq 2.$So we can consider $\mu_{m-b}.$ 
Note that $2\in B$ and so $B\neq \emptyset.$  Let $k''=\max(B),$ since $\mu_{m-1}+\mu_{m-{\mu_{m-1}}}\geq m$ we have $k''\neq \mu_{m-1}$ and so  $k''+1\leq \mu_{m-1}.$ Since $k''=\max(B),$ on the one hand we have that $k''+1\notin B$ and $m+1\leq \mu_{m-k''-1}+(k''+1)+1.$ On the other hand $\mu_{m-k''}+k''+1<m+1.$ So $k=k''+1$ is the integer we seek. 
\\
To prove that this integer is unique, let us assume there exist two positive integers $2\leq k_1<k_2\leq \mu_{m-1}$ such that $\mu$ transforms to: $$\Lambda_i:(\mu_1+1,\cdots,\mu_{m-k_i}+1,\mu_{m-k_i+1},\cdots,\mu_{m-1},k_i),$$ 
such that $k_i+\mu_{m-k_i+1}<m+1\leq k_i+1+\mu_{m-k_i},$ for $i=1,2.$ Since $k_1<k_2,$ we have $m-k_1\geq m-k_2+1$ and so $\mu_{m-k_1}\leq \mu_{m-k_2+1}.$ So we have $k_1+\mu_{m-k_1}+1\leq k_2+\mu_{m-k_2+1}<m+1,$ which is a contradiction and proves the last bijection.
\\
\\
These three bijections give us a one-to-one correspondence between the partitions counted by $c_{3,2}(m,n)-c_{3,1}(m,n)$ and those counted by $c_{3,2}(m-1,n-m).$ 
 \\
 \\
In order to prove the last equation of our system let $\Lambda:(\lambda_1,\cdots,\lambda_m)$ be a partition of $n$ such that $\lambda_m \neq 1$ and either $\lambda_m\geq m+1$ or $\lambda_m \leq m$ and $\lambda_m+\lambda_{m+1-\lambda_m}\geq m+2.$ We transform this partition to another partition $\mu,$ by subtracting $1$ from each part. So $\mu_m=\lambda_{m}-1\geq 1.$ If $\lambda_m\geq m+1$ then we have $\mu_m =\lambda_m-1\geq m$. So if $\mu_m=1$ then $\mu:(1).$
\\
If $\lambda_m <m+1$ %on the one hand  
we have: $$\mu_m+\mu_{m-\mu_m}=\lambda_m-1+\lambda_{m-\lambda_m+1}-1\geq (m+2)-2=m.$$
%On the other hand, if $\mu_m=\mu_{m-1}=\mu_{m-2}=1,$ then $\lambda_m=\lambda_{m-1}=\lambda_{m-2}=2.$ So we have $m+2\leq \lambda_m+\lambda_{m+1-\lambda_m}=2+\lambda_{m-1}=2+2=4.$ This implies $m\leq 2,$ which is a contradiction because we supposed that $\Lambda$ has at least $3$ parts: $\lambda_m,\lambda_{m-1}$ and $\lambda_{m-2}.$ this means that $\mu$ can have at most two parts equal to $1.$ 
So we have defined a transformation from $\Lambda$ to $\mu$ which is clearly a bijection. 
\\
\\
Now let $b_{3,i}(m,n)$ denote the number of partitions of $n$ with exactly $m$ parts and which are counted by $B_{3,i}(n)$.Then the $b_{3,i}(m,n)$ are \textit {uniquely} determined by Andrew's system of equations (see the proof of Theorem \ref {Gordon}). Therefore, $c_{3,i} (m,n)=b_{3,i} (m,n)$ for all $m$ and $n$ with $0 \leq i \leq 3.$
\\
\\
Since $\sum_{m \geq 0} c_{3,i} (m,n)= C_{3,i}(n)$ and $\sum_{m \geq 0} b_{3,i} (m,n)= B_{3,i}(n)$, we have

$$C_{3,i}(n)=\sum_{m \geq 0} c_{3,i} (m,n)=\sum_{m \geq 0} b_{3,i} (m,n)= B_{3,i}(n)=A_{3,i}(n).$$
\end{proof}

\section{ANALYTIC FORM OF A NEW VERSION OF GORDON'S IDENTITIES FOR THE CASE  $r=3$}

In this section we give an analytic form of Theorem \ref{r3} for $r=i=3$. Let $m,c\geq 1$ and $r\geq2$.

  We define $r$ blocks of integers which are greater than or equal to $m$ as follows: The first block contains $c$ integers. The number of integers which appear in each block is equal to the last number of the previous block. i.e.,
 
 $$ \underbrace{ n_{1,1} \leq \cdots \leq n_{1,c} }_{\text {The first block}}\leq \underbrace{ n_{2,1} \leq \cdots \leq n_{2,n_{1,c}} }_{\text {The second block}}\leq \underbrace{ n_{3,1} \leq \cdots \leq n_{3,n_{2,n_{1,c}}}}_{\text{The third block}}\leq \cdots.$$
 
In order to simplify notations, for $1\leq j \leq r$ we introduce:

$$f(j)= \begin{cases}
c &\text{ if } j=1 \\
n_{j-1,f(j-1)} &\text{ if } j\geq 2.

\end{cases}$$

So we are considering the following $r$ blocks of positive integers:

  $$\underbrace{ n_{1,1} \leq \cdots \leq n_{1,f(1)} }_{\text {The first block}}\leq \underbrace{ n_{2,1} \leq \cdots \leq n_{2,f(2)} }_{\text {The second block}}\leq \underbrace{ n_{3,1} \leq \cdots \leq n_{3,f(3)}}_{\text{The third block}}\leq \cdots \leq \underbrace{ n_{r,1} \leq \cdots \leq n_{r,f(r)} }_{\text {The $r$-th block}}.$$

 We denote by $H_{r,c}^m$ the Hilbert-Poincar\'e series of the following algebra:

$$\frac{\mathbf k[x_m,x_{m+1},\cdots]}{(x_{n_{1,1}}\cdots x_{n_{1,f(1)}}x_{n_{2,1}}\cdots x_{n_{2,f(2)}}\cdots x_{n_{r,1}} \cdots x_{n_{r,f(r)}})}.$$

\begin{lem}\label{H_{2,c^m}} We have:
$$H_{2,c}^m=\sum_{n=0}^{m-1}\frac{q^{nm}}{(q)_n}+\sum_{j=0}^{c-1}\sum_{m\leq \ell_j \leq \cdots \leq \ell_1\leq k}\frac{q^{k^2+\ell_j+\cdots+\ell_1}}{(q)_k}.$$
%$$H_{2,c}^m=\sum_{n=0}^{m-1}\frac{q^{nm}}{(q)_n}+\sum_{n\geq m}\frac{q^{n^2}}{(q)_n}+\sum_{j=1}^{c-1}\sum_{m\leq \ell_j \leq \cdots \leq \ell_1\leq n}\frac{q^{n^2+\ell_j+\cdots+\ell_1}}{(q)_n}.$$
\end{lem}
% braye $H_{2,c^m}$: 2= tedade blocks, c=tedade azaye avalin block, m ham ke vazehe.
\begin{proof}
The proof is by induction on $c.$ We obtain the proof for the case $c=1$ by the same computation as in Theorem $1.5$ from \cite{AM}, considering $H_{2,1}^m$ instead of $H_{2,1}^1$ (note that in \cite{AM} we denoted $H_{2,1}^m$ by $H_{m}$). Suppose that the equality holds for $H_{2,c}^l$ for $l\geq 1$. For $c+1$ we have:

  $$H_{2,c+1}^m=1+\sum_{l \geq m}q^l H_{2,c}^{l},$$
 which, by the induction hypothesis is equal to:
 
\begin{equation*}
\begin{aligned}
1+\sum_{l\geq m} q^l\Big(\sum_{n=0}^{l-1}\frac{q^{nl}}{(q)_n}+\sum_{j=0}^{c-1}\sum_{l\leq \ell_j \leq \cdots \leq \ell_1\leq k}\frac{q^{k^2+\ell_j+\cdots+\ell_1}}{(q)_k}\Big)\\
=1+\sum_{l\geq m } \sum_{n=0}^{l-1}\frac{q^{(n+1)l}}{(q)_n}+\sum_{l \geq m} \sum_{j=0}^{c-1}\sum_{l\leq \ell_j \leq \cdots \leq \ell_1\leq k}\frac{q^{k^2+\ell_j+\cdots+\ell_1+l}}{(q)_k}.
\end{aligned}
\end{equation*}

  %$$q^m\Big(1+\sum_{n=1}^{m-1}\frac{q^{nm}}{(q)_n}+\sum_{n\geq m}\frac{q^{n^2}}{(q)_n}+\sum_{j=1}^{c-1}\sum_{m\leq \ell_j \leq \cdots \leq \ell_1\leq n}\frac{q^{n^2+\ell_1+\cdots+\ell_j}}{(q)_n}\Big)+$$
  
 % $$q^{m+1}\Big(1+\sum_{n=1}^{m}\frac{q^{n(m+1)}}{(q)_n}+\sum_{n\geq m+1}\frac{q^{n^2}}{(q)_n}+\sum_{j=1}^{c-1}\sum_{m+1\leq \ell_j \leq \cdots \leq \ell_1\leq n}\frac{q^{n^2+\ell_1+\cdots+\ell_j}}{(q)_n}\Big)+$$
  
  %$$q^{m+2}\Big(1+\sum_{n=1}^{m+1}\frac{q^{n(m+2)}}{(q)_n}+\sum_{n\geq m+2}\frac{q^{n^2}}{(q)_n}+\sum_{j=1}^{c-1}\sum_{m+2\leq \ell_j \leq \cdots \leq \ell_1\leq n}\frac{q^{n^2+\ell_1+\cdots+\ell_j}}{(q)_n}\Big)+$$
  
  %$$\vdots$$
 % For $1\leq i \leq 4$ we denote by $\mathcal{C}_i$ the sum of the $i$-th column in the  series above.
% Obviously the first term is equal to $\frac{q^m}{1-q}.$
 %The sum of the second column $\mathcal{C}_2,$ is equal to:
% $$\sum_{n=1}^{m-1}\frac{q^{(n+1)m}}{(q)_n}+\sum_{n=1}^{m}\frac{q^{(n+1)(m+1)}}{(q)_n}+\sum_{n=1}^{m+1}\frac{q^{(n+1)(m+2)}}{(q)_n}+\cdots$$

  %Which can write as:
  
%$$\begin{matrix}
%  \frac{q^{2m}}{(q)_1} & +\frac{q^{3m}}{(q)_2} & +\cdots  & +\frac{q^{m^2-m}}{(q)_{m-2}}&+\frac{q^{m^2}}{(q)_{m-1}}&&& \\
   %+  \frac{q^{2m+2}}{(q)_1}  &+ \frac{q^{3m+3}}{(q)_2}  & +\cdots & +\frac{q^{m^2-1}}{(q)_{m-2}}& +\frac{q^{m^2+m}}{(q)_{m-1}}&+\frac{q^{(m+1)^2}}{(q)_{m}}+&& \\
%    + \frac{q^{2m+4}}{(q)_1}&+\frac{q^{3m+6}}{(q)_2}&+\cdots &+\frac{q^{m^2+m-2}}{(q)_{m-2}} &+\frac{q^{m^2+2m}}{(q)_{m-1}}&+\frac{q^{m^2+3m+2}}{(q)_{m}}&+\frac{q^{(m+2)^2}}{(q)_{m+1}}& \\
% +~~~\cdots\\   
% \end{matrix}$$

  %Summing by columns we obtain that this sum is equal to
 By inverting the summation in the second term we obtain:
  %$\frac{q^{2m}}{(q)_2}+\frac{q^{3m}}{(q)_3}+\cdots+\frac{q^{m^2-m}}{(q)_{m-1}}+\frac{q^{m^2}}{(q)_m}+\frac{q^{(m+1)^2}}{(q)_{m+1}}+\frac{q^{(m+2)^2}}{(q)_{m+2}}+\cdots$
  
  %  $=\sum_{n=2}^{m-1}\frac{q^{nm}}{(q)_n}+\sum_{n\geq m } \frac{q^{n^2}}{(q)_n}.$
$$\sum_{l\geq m } \sum_{n=0}^{l-1}\frac{q^{(n+1)l}}{(q)_n}
=\sum_{n=1}^{m-1}\frac{q^{nm}}{(q)_n}+\sum_{n\geq m } \frac{q^{n^2}}{(q)_n}.$$

% We now compute $\mathcal{C}_3$:
 
% $\mathcal{C}_3=\sum_{n \geq m}\frac{q^{n^2+m}}{(q)_n}+\sum_{n\geq m+1}\frac{q^{n^2+m+1}}{(q)_n}+\sum_{n\geq m+2}\frac{q^{n^2+m+2}}{(q)_n}+\cdots$
 
% $=\sum_{m\leq \ell \leq n}\frac{q^{n^2+\ell}}{(q)_n}.$
%The sum of the last column $\mathcal{C}_4$ is equal to:
Changing $l$ by $\ell_{j+1}$ and then $j+1$ by $j$ in the last term we obtain:
% $\sum_{j=1}^{c-1}\sum_{m\leq \ell_j \leq \cdots \leq \ell_1\leq n}\frac{q^{n^2+m+\ell_1+\cdots+\ell_j}}{(q)_n}+\sum_{j=1}^{c-1}\sum_{m+1\leq \ell_j \leq \cdots \leq \ell_1\leq n}\frac{q^{n^2+(m+1)+\ell_1+\cdots+\ell_j}}{(q)_n}$
 
% $+\sum_{j=1}^{c-1}\sum_{m+2\leq \ell_j \leq \cdots \leq \ell_1\leq n}\frac{q^{n^2+(m+2)+\ell_1+\cdots+\ell_j}}{(q)_n}+\cdots$

\begin{equation*}
\begin{aligned}
\sum_{l \geq m} \sum_{j=0}^{c-1}\sum_{l\leq \ell_j \leq \cdots \leq \ell_1\leq k}\frac{q^{k^2+\ell_j+\cdots+\ell_1+l}}{(q)_k}&
=\sum_{j=0}^{c-1}\sum_{m\leq \ell_{j+1}\leq \ell_j \leq \cdots \leq \ell_1\leq k}\frac{q^{k^2+\ell_1+\cdots+\ell_{j+1}}}{(q)_k}\\
&=\sum_{j=1}^{c}\sum_{m\leq \ell_j \leq \cdots \leq \ell_1\leq k}\frac{q^{k^2+\ell_1+\cdots+\ell_j}}{(q)_k}.
\end{aligned}
\end{equation*}
% $\sum_{j=1}^{c-1}\sum_{m\leq \ell_{j+1\leq \ell_j} \leq \cdots \leq \ell_1\leq n}\frac{q^{n^2+\ell_1+\cdots+\ell_{j+1}}}{(q)_n}\underset{\text{changing } j+1\ \text{by}\ j }{=}\sum_{j=2}^{c}\sum_{m\leq \ell_j \leq \cdots \leq \ell_1\leq n}\frac{q^{n^2+\ell_1+\cdots+\ell_j}}{(q)_n}.$

% So far we have computed $\mathcal{C}_i$ for $1\leq i \leq 4.$ So we have:
 So $H_{2,c+1}^m$ is equal to:
 
% $$H_{2,c+1}^m=1+\mathcal{C}_1+\mathcal{C}_2+\mathcal{C}_3+\mathcal{C}_4=$$
 
% $$1+(\frac{q^m}{(q)_1})+\Big(\sum_{n=2}^{m-1}\frac{q^{nm}}{(q)_n}+\sum_{n\geq m } \frac{q^{n^2}}{(q)_n}\Big)+(\sum_{m\leq \ell \leq n}\frac{q^{n^2+\ell}}{(q)_n})+\Big(\sum_{j=2}^{c}\sum_{m\leq \ell_j \leq \cdots \leq \ell_1\leq n}\frac{q^{n^2+\ell_1+\cdots+\ell_j}}{(q)_n}\Big)=$$
 
% $$1+\sum_{n=1}^{m-1}\frac{q^{nm}}{(q)_n}+\sum_{n\geq m}\frac{q^{n^2}}{(q)_n}+\sum_{j=1}^{c}\sum_{m\leq \ell_j \leq \cdots \leq \ell_1\leq n}\frac{q^{n^2+\ell_1+\cdots+\ell_j}}{(q)_n}.$$

 \begin{equation*}
\begin{aligned}
1+\Big(\sum_{n=1}^{m-1}\frac{q^{nm}}{(q)_n}+\sum_{n\geq m } \frac{q^{n^2}}{(q)_n}\Big)+\Big(\sum_{j=1}^{c}\sum_{m\leq \ell_j \leq \cdots \leq \ell_1\leq k}\frac{q^{k^2+\ell_1+\cdots+\ell_j}}{(q)_k}\Big)\\
=\sum_{n=0}^{m-1}\frac{q^{nm}}{(q)_n}+\sum_{j=0}^{c}\sum_{m\leq \ell_j \leq \cdots \leq \ell_1\leq k}\frac{q^{k^2+\ell_1+\cdots+\ell_j}}{(q)_k}.
\end{aligned}
\end{equation*}
\end{proof}

\begin{prop}\label{HP3} We have: 

$$HP_{\frac{S}{I'_{3,3}}}(q)=\sum_{0\leq j \leq \ell_j \leq \cdots \leq \ell_1 \leq n} \frac{q^{n^2+\ell_1+\cdots+\ell_j}}{(q)_n},$$
where $$I'_{3,3}=( x_c x_{k_1}\cdots x_{k_c} x_{i_1}\cdots x_{i_{k_c}}| 1\leq c\leq k_1 \leq \cdots \leq k_c \leq i_{k_c} \leq \cdots \leq i_1).$$
%where $I'_{3,3}$ is the ideal introduced in the sixth section. 
\end{prop}

\begin{proof}
We denote as usual $HP_{\frac{S}{I'_{3,3}}}(q)$ by $H_{3,1}^1.$
Using repetitively Equation (\ref{eq:1}) we obtain:
%$$H_{3,1}^1=1+qH_{2,1}^1+q^2H_{2,2}^2+q^3H_{2,3}^3+q^4H_{2,4}^4+\cdots$$

$$H_{3,1}^1=1+\sum_{l\geq 1}q^l H^l_{2,l},$$
which gives by Lemma \ref{H_{2,c^m}}, inverting summations, shifting indices and easy computations:

\begin{equation*}
\begin{aligned}
H_{3,1}^1&=1+\sum_{l \geq 1}\sum_{n=0}^{l-1}\frac{q^{(n+1)l}}{(q)_n}+\sum_{l\geq 1} \sum_{j=0}^{l-1}\sum_{l\leq \ell_j \leq \cdots \leq \ell_1\leq k}\frac{q^{k^2+\ell_j+\cdots+\ell_1+l}}{(q)_k}\\
&=1+\sum_{n\geq 1}\frac{q^{n^2}}{(q)_n}+\sum_{j=1}^{\ell_j}\sum_{1\leq \ell_j \leq \cdots \leq \ell_1\leq k}\frac{q^{k^2+\ell_j+\cdots+\ell_1}}{(q)_k}\\
&=\sum_{0\leq j \leq \ell_j \leq \cdots \leq \ell_1 \leq n} \frac{q^{n^2+\ell_1+\cdots+\ell_j}}{(q)_n}.
\end{aligned}
\end{equation*}

 %\begin{equation*}
%\begin{aligned}
%H_{3,1}^1=1+(\frac{q}{1-q})+(\sum_{n\geq 2} \frac{q^{n^2}}{(q)_n})+\Big(\sum_{1\leq \ell \leq n}\frac{q^{n^2+\ell}}{(q)_n}+\sum_{2\leq j\leq \ell_j \leq \cdots \leq \ell_1\leq n}\frac{q^{n^2+\ell_1+\cdots+\ell_j}}{(q)_n}\Big)\\
%=\sum_{n\geq 0}\frac{q^{n^2}}{(q)_n}+\sum_{1\leq j \leq \ell_j \leq \cdots \leq \ell_1 \leq n} \frac{q^{n^2+\ell_1+\cdots+\ell_j}}{(q)_n}.
%\end{aligned}
%\end{equation*}
\end{proof}

\begin{theo}\label{anal 3} We have 
$$\sum_{0\leq n_1, n_2} \frac{q^{(n_1+n_2)^2+n_2^2}}{(q)_{n_1}(q)_{n_2}}=\sum_{0\leq j \leq \ell_j \leq \cdots \leq \ell_1 \leq n} \frac{q^{n^2+\ell_1+\cdots+\ell_j}}{(q)_n}.$$
\end{theo}

\begin{proof}
By Proposition \ref{HP3} the right hand side of this equation is equal to $H_{3,1}^1$ which is the Hilbert-Poincar\'e series of the graded algebra $\frac{\mathbf k[x_1,x_2,\cdots]}{I'_{3,3}}.$ By Proposition \ref{I'}, this Hilbert-Poincar\'e series is the generating series for the partitions counted by $C_{3,3}(n).$ 
\\
Note also that the left hand side of the above equality is the generating series of the partitions  counted by $B_{3,3}(n)$ in Gordon's identities (see Theorem \ref{cor}). 
\\
But by Theorem \ref{r3} we know that $B_{3,3}(n)=C_{3,3}(n).$ So we have:

%$$\sum_{0\leq n_1, n_2} \frac{q^{(n_1+n_2)^2+n_2^2}}{(q)_{n_1}(q)_{n_2}}=\sum_{n \geq 0} B_{3,3}(n)q^n=$$
%$$\sum_{n\geq 0} C_{3,3}(n)q^n=1+\sum_{n\geq 1}\frac{q^{n^2}}{(q)_n}+\sum_{1\leq j \leq \ell_j \leq \cdots \leq \ell_1 \leq n} \frac{q^{n^2+\ell_1+\cdots+\ell_j}}{(q)_n}.$$
\begin{equation*}
\begin{aligned}
\sum_{0\leq n_1, n_2} \frac{q^{(n_1+n_2)^2+n_2^2}}{(q)_{n_1}(q)_{n_2}} & =\sum_{n \geq 0} B_{3,3}(n)q^n \\
      &=\sum_{n\geq 0} C_{3,3}(n)q^n \\
      &=\sum_{0\leq j \leq \ell_j \leq \cdots \leq \ell_1 \leq n} \frac{q^{n^2+\ell_1+\cdots+\ell_j}}{(q)_n}.
\end{aligned}
\end{equation*}

\end{proof}

\begin{rem}
By this theorem and the analytic form of Gordon's identities (see Theorem \ref{cor}) we obtain that the series which appear in this theorem are also equal to $\underset{\underset{n \not\equiv 0, \pm 3 (mod. 7)}{n\geq1}}{\prod} \frac{1}{1-q^n},$ which is the generating series of the partitions counted by $A_{3,3}(n)$ in Gordon's identities.
\end{rem}
Recall that the \textit{q-binomial numbers} ${\binom{ N+m}{ m}}_q$ are defined as the generating series of the integer partitions with length $\leq m$ and each part $\leq N,$ which is equal to:
$${\binom{ N+m}{ m}}_q=\frac{(q)_{N+m}}{(q)_m (q)_{N}},$$
where $(q)_j=(1-q)\cdots(1-q^j).$
We now give a direct proof of Theorem \ref{anal 3}. To do so we need the following lemma:

\begin{lem}\label{lem} For all integers $n\geq 1$ and $0\leq j \leq n$ we have: 
$$ \sum_{j\leq \ell_j \leq \cdots \leq \ell_1 \leq n} q^{\ell_1+\cdots+\ell_j-j^2}={\binom{n}{j}}_q.$$
%$$(q)_j \sum_{j\leq \ell_j \leq \cdots \leq \ell_1 \leq n} q^{\ell_1+\cdots+\ell_j-j^2}=(1-q^{n-j+1})(1-q^{n-j+2})\cdots (1-q^n);$$

\end{lem}

\begin{proof} 

Shifting all the indices $\ell_i$ to $\ell_i-j$ on the left-hand side of the above equation gives:
$$ \sum_{0\leq \ell_j \leq \cdots \leq \ell_1 \leq n-j} q^{\ell_1+\cdots+\ell_j};$$ 
which is the generating series of integer partitions of length $\leq j$ and each part $\leq n-j$ and so by definition it is equal to ${\binom{n}{j}}_q.$

\end{proof}

\begin{theo}\label{direct}
We have:

$$\sum_{ n_1, n_2\geq 0} \frac{q^{(n_1+n_2)^2+n_2^2}}{(q)_{n_1}(q)_{n_2}}=\sum_{0 \leq j \leq \ell_j \leq \cdots \leq \ell_1 \leq n} \frac{q^{n^2+\ell_1+\cdots+\ell_j}}{(q)_n}.$$
\end{theo}

\begin{proof} Changing $n_2$ by $j$ and $n_1+n_2$ by $n$ in the left hand side of the equation above we obtain: 
\begin{equation*}
\begin{aligned}
\sum_{ n_1, n_2\geq 0} \frac{q^{(n_1+n_2)^2+n_2^2}}{(q)_{n_1}(q)_{n_2}}& = \sum_{n\geq j \geq 0} \frac{q^{n^2+j^2}}{(q)_{n-j}(q)_j} \\
      &=\sum_{n\geq j \geq 0} \frac{q^{n^2+j^2}}{(q)_n}{\binom{n}{j}}_q\\
      &= \sum_{n\geq 0}\frac{q^{n^2}}{(q)_n}+\sum_{n\geq j \geq 1} \frac{q^{n^2+j^2}}{(q)_n}{\binom{n}{j}}_q.
\end{aligned}
\end{equation*}

By Lemma \ref{lem} this is equal to:

$$
\sum_{n\geq 0}\frac{q^{n^2}}{(q)_n}+\sum_{n\geq j \geq 1} \frac{q^{n^2+j^2}}{(q)_n}  \sum_{j\leq \ell_j \leq \cdots \leq \ell_1 \leq n} q^{\ell_1+\cdots+\ell_j-j^2} \\
      =\sum_{0 \leq j \leq \ell_j \leq \cdots \leq \ell_1 \leq n} \frac{q^{n^2+\ell_1+\cdots+\ell_j}}{(q)_n}.$$

\end{proof}

\section{ANALYTIC FORM OF THE CONJECTURE }
In this section we keep the same notations. We take $i=r$ all along this section. In order to give an analytic form of Conjecture \ref{conjG} in this case, we need the following lemma:

\begin{lem}\label{H_{3,c}^m} We have:
$$H_{3,c}^m=\sum_{n=0}^{m-1}\frac{q^{nm}}{(q)_n}+\sum_{n\geq m}\frac{q^{n^2}}{(q)_n}+\sum_{\underset{1\leq j \leq \ell_{j-c+1}+c-1}{m\leq \ell_j \leq \cdots \leq \ell_1 \leq n}} \frac{q^{n^2+\ell_1+\cdots +\ell_j}}{(q)_n}.$$
\end{lem}

\begin{proof}
The proof is by induction on $c.$ In Proposition \ref{HP3} we proved the case $c=1$ of this induction. Suppose that the equality holds for $H_{3,c}^l$ for all $l\geq 1$ and we prove it for $H_{3,c+1}^m.$ Using repetitively Equation (\ref{eq:1}) we have:
%$$H_{3,c+1}^m=1+q^mH_{3,c}^m+q^{m+1}H_{3,c}^{m+1}+q^{m+2}H_{3,c}^{m+2}+\cdots$$
$$H_{3,c+1}^m=1+\sum_{l\geq m} q^l H_{3,c}^{l},$$
 which by the induction hypothesis is equal to:

%$$\begin{matrix}

  %    1 & +q^m\Big(1   & +\sum_{n=1}^{m-1}\frac{q^{nm}}{(q)_n} & +\sum_{n\geq m}\frac{q^{n^2}}{(q)_n}& +\sum_{\underset{1\leq j \leq \ell_{j-c+1}+c-1}{m\leq \ell_j \leq \cdots \leq \ell_1 \leq n}} \frac{q^{n^2+\ell_1+\cdots +\ell_j}}{(q)_n} \Big)& \\
 %      & +q^{m+1}\Big(1  & +\sum_{n=1}^{m}\frac{q^{n(m+1)}}{(q)_n} & +\sum_{n\geq m+1}\frac{q^{n^2}}{(q)_n}& +\sum_{\underset{1\leq j \leq \ell_{j-c+1}+c-1}{m+1\leq \ell_j \leq \cdots \leq \ell_1 \leq n}} \frac{q^{n^2+\ell_1+\cdots +\ell_j}}{(q)_n} \Big)& \\
  %    & +q^{m+2}\Big(1  & +\sum_{n=1}^{m+1}\frac{q^{n(m+2)}}{(q)_n} & +\sum_{n\geq m+2}\frac{q^{n^2}}{(q)_n}& +\sum_{\underset{1\leq j \leq \ell_{j-c+1}+c-1}{m+2\leq \ell_j \leq \cdots \leq \ell_1 \leq n}} \frac{q^{n^2+\ell_1+\cdots +\ell_j}}{(q)_n} \Big)&\\            
  %    &  +~~~\cdots\\
  %     \end{matrix}$$
  \begin{equation*}
\begin{aligned}
1+\sum_{l\geq m} q^l\Big(\sum_{n=0}^{l-1}\frac{q^{nl}}{(q)_n}+\sum_{n\geq l}\frac{q^{n^2}}{(q)_n}+\sum_{\underset{1\leq j \leq \ell_{j-c+1}+c-1}{l\leq \ell_j \leq \cdots \leq \ell_1 \leq n}} \frac{q^{n^2+\ell_1+\cdots +\ell_j}}{(q)_n}\Big)\\
=1+\sum_{l\geq m } \sum_{n=0}^{l-1}\frac{q^{(n+1)l}}{(q)_n}+\sum_{l \geq m}\sum_{n\geq l}\frac{q^{n^2+l}}{(q)_n}+\sum_{\underset{1\leq j \leq \ell_{j-c+1}+c-1}{l\leq \ell_j \leq \cdots \leq \ell_1 \leq n}} \frac{q^{n^2+\ell_1+\cdots +\ell_j+l}}{(q)_n}.
\end{aligned}
\end{equation*}       
       
    %  For $1\leq i \leq 4$ we denote by $\mathcal{C}_i$ the sum of the $i$th column of the  series above. 
    Doing similar computations as we did in the proof of Lemma \ref{H_{2,c^m}} we obtain that $H_{3,c}^m$ is equal to:

% $\mathcal{C}_1=\frac{q^m}{(q)_1};\ \mathcal{C}_2=\sum_{n=2}^{m-1}\frac{q^{nm}}{(q)_n}+\sum_{n\geq m } \frac{q^{n^2}}{(q)_n};\ \mathcal{C}_3=\sum_{m\leq \ell \leq n}\frac{q^{n^2+\ell}}{(q)_n};\ \mathcal{C}_4= \sum_{\underset{2\leq s \leq \ell_{s-c}+c}{m\leq \ell_{s} \leq \cdots \leq \ell_1 \leq n}} \frac{q^{n^2+\ell_1+\cdots +\ell_s}}{(q)_n}.$
  
% So we have $H_{3,c}^m=1+\mathcal{C}_1+\mathcal{C}_2+\mathcal{C}_3+\mathcal{C}_4$ which is equal to:
 
$$1+\Big(\sum_{n=1}^{m-1}\frac{q^{nm}}{(q)_n}+\sum_{n\geq m } \frac{q^{n^2}}{(q)_n}\Big)+(\sum_{m\leq \ell \leq n}\frac{q^{n^2+\ell}}{(q)_n})+\Big(\sum_{\underset{2\leq j \leq \ell_{j-c}+c}{m\leq \ell_{j} \leq \cdots \leq \ell_1 \leq n}}\frac{q^{n^2+\ell_1+\cdots +\ell_j}}{(q)_n}\Big)$$
 
 $$=\sum_{n=0}^{m-1}\frac{q^{nm}}{(q)_n}+\sum_{n\geq m}\frac{q^{n^2}}{(q)_n}+\sum_{\underset{1\leq j \leq \ell_{j-c+1}+c-1}{m\leq \ell_j \leq \cdots \leq \ell_1 \leq n}} \frac{q^{n^2+\ell_1+\cdots +\ell_j}}{(q)_n}.$$

 \end{proof}

 \begin{prop} Given integers $m,c\geq 1$ and $r\geq 3$ we have:
 
 $$H_{r,c}^m=\sum_{n=0}^{m-1}\frac{q^{nm}}{(q)_n}+\sum_{n\geq m}\frac{q^{n^2}}{(q)_n}+\sum_{\underset{1\leq j \leq \sum_{i=1}^{r-2}p_{r,i}(\mu)+c-1}{\underset{\mu=\ell_1+\cdots+\ell_{j-c+1}}{m\leq \ell_j \leq \cdots \leq \ell_1 \leq n}}} \frac{q^{n^2+\ell_1+\cdots +\ell_j}}{(q)_n},$$
 where $p_{r,i}(\mu)$ is the $(r,i)$-new part of $\mu$ (see Definition \ref{def}).
  \end{prop}

\begin{proof}
The proof is by induction on $r.$ By Lemma \ref{H_{3,c}^m} the case $r=3$ of the induction is true. Assume that the equality holds for $H_{r,c}^m.$ We prove it for $H_{r+1,c}^m$ using induction on $c$. Using repetitively Equation (\ref{eq:1}) for $c=1$ we have:

%$$H_{r+1,1}^m=1+q^mH_{r,m}^m+q^{m+1}H_{r,m+1}^{m+1}+q^{m+2}H_{r,m+2}^{m+2}+\cdots$$

$$H_{r+1,1}^m=1+\sum_{l\geq m}q^lH^l_{r,l}.$$
By the induction hypothesis on $r$ it is equal to:

%$$\begin{matrix}

  %    1  +q^m\Big(1   & +\sum_{n=1}^{m-1}\frac{q^{nm}}{(q)_n} & +\sum_{n\geq m}\frac{q^{n^2}}{(q)_n}& +\sum_{\underset{1\leq j \leq \sum_{i=1}^{r-2}p_{r,i}(\mu)+m-1}{\underset{\mu=\ell_1+\cdots+\ell_{j-m+1}}{m\leq \ell_j \leq \cdots \leq \ell_1 \leq n}}} \frac{q^{n^2+\ell_1+\cdots +\ell_j}}{(q)_n} \Big)& \\
 %       +q^{m+1}\Big(1   & +\sum_{n=1}^{m}\frac{q^{n(m+1)}}{(q)_n} & +\sum_{n\geq m+1}\frac{q^{n^2}}{(q)_n}& +\sum_{\underset{1\leq j \leq \sum_{i=1}^{r-2}p_{r,i}(\mu)+m}{\underset{\mu=\ell_1+\cdots+\ell_{j-m+2}}{m+1\leq \ell_j \leq \cdots \leq \ell_1 \leq n}}} \frac{q^{n^2+\ell_1+\cdots +\ell_j}}{(q)_n} \Big)& \\
   %     +q^{m+2}\Big(1   & +\sum_{n=1}^{m+1}\frac{q^{n(m+2)}}{(q)_n} & +\sum_{n\geq m+2}\frac{q^{n^2}}{(q)_n}& +\sum_{\underset{1\leq j \leq \sum_{i=1}^{r-2}p_{r,i}(\mu)+m+1}{\underset{\mu=\ell_1+\cdots+\ell_{j-m+3}}{m+2\leq \ell_j \leq \cdots \leq \ell_1 \leq n}}} \frac{q^{n^2+\ell_1+\cdots +\ell_j}}{(q)_n} \Big)&\\            
  %      +~~~\cdots\\ 
  %      \end{matrix}$$
        
          \begin{equation}\label{in}
\begin{aligned}
1+\sum_{l\geq m} q^l\Big(\sum_{n=0}^{l-1}\frac{q^{nl}}{(q)_n}+\sum_{n\geq l}\frac{q^{n^2}}{(q)_n}+\sum_{\underset{1\leq j \leq \sum_{i=1}^{r-2}p_{r,i}(\mu)+l-1}{\underset{\mu=\ell_1+\cdots+\ell_{j-l+1}}{l\leq \ell_j \leq \cdots \leq \ell_1 \leq n}}} \frac{q^{n^2+\ell_1+\cdots +\ell_j}}{(q)_n}\Big)\\
=1+\sum_{l\geq m } \sum_{n=0}^{l-1}\frac{q^{(n+1)l}}{(q)_n}+\sum_{l \geq m}\sum_{n\geq l}\frac{q^{n^2+l}}{(q)_n}+\sum_{\underset{1\leq j \leq \sum_{i=1}^{r-2}p_{r,i}(\mu)+l-1}{\underset{\mu=\ell_1+\cdots+\ell_{j-l+1}}{l\leq \ell_j \leq \cdots \leq \ell_1 \leq n}}} \frac{q^{n^2+\ell_1+\cdots +\ell_j+l}}{(q)_n}.
\end{aligned}
\end{equation}

%If we denote by $\mathcal{C}_i$ the sum of the $i$-th column of the  equation above, by the proof of Lemma \ref{H_{2,c^m}} we know that:
%$\mathcal{C}_1=\frac{q^m}{(q)_1};$
%$\mathcal{C}_2=\sum_{n=2}^{m-1}\frac{q^{nm}}{(q)_n}+\sum_{n\geq m } \frac{q^{n^2}}{(q)_n};$
%$\mathcal{C}_3=\sum_{m\leq \ell \leq n}\frac{q^{n^2+\ell}}{(q)_n}.$
%We have:

By the proof of Lemma \ref{H_{2,c^m}} we know that the sum of the first three terms of Equation \ref{in}  is equal to:

$$\sum_{n=1}^{m-1}\frac{q^{nm}}{(q)_n}+\sum_{n\geq m } \frac{q^{n^2}}{(q)_n}+\sum_{m\leq \ell \leq n}\frac{q^{n^2+\ell}}{(q)_n}.$$

Changing $l$ by $\ell_{j+1}$ and then $j+1$ by $s$ in the last term of this equation gives us:
%$ \mathcal{C}_4= \sum_{\underset{1\leq j \leq \sum_{i=1}^{r-2}p_{r,i}(\mu)+m-1}{\underset{\mu=\ell_1+\cdots+\ell_{j-m+1}}{m\leq \ell_j \leq \cdots \leq \ell_1 \leq n}}} \frac{q^{n^2+\ell_1+\cdots +\ell_j+m}}{(q)_n}+ \sum_{\underset{1\leq j \leq \sum_{i=1}^{r-2}p_{r,i}(\mu)+m}{\underset{\mu=\ell_1+\cdots+\ell_{j-m+2}}{m+1\leq \ell_j \leq \cdots \leq \ell_1 \leq n}}} \frac{q^{n^2+\ell_1+\cdots +\ell_j+(m+1)}}{(q)_n}$

%$$+\sum_{\underset{1\leq j \leq \sum_{i=1}^{r-2}p_{r,i}(\mu)+m+1}{\underset{\mu=\ell_1+\cdots+\ell_{j-m+3}}{m+2\leq \ell_j \leq \cdots \leq \ell_1 \leq n}}} \frac{q^{n^2+\ell_1+\cdots +\ell_j+(m+2)}}{(q)_n}+ \cdots $$

%Which is equal to:

%$\sum_{\underset{1\leq j \leq \sum_{i=1}^{r-2}p_{r,i}(\mu)+\ell_{j+1}-1}{\underset{\mu=\ell_1+\cdots+\ell_{j-\ell_{j+1}+1}}{m\leq \ell_{j+1} \leq \ell_j \leq \cdots \leq \ell_1 \leq n}}} \frac{q^{n^2+\ell_1+\cdots +\ell_j+\ell_{j+1}}}{(q)_n}\underset{\text{changing } j+1\ \text{by}\ s }{=}\sum_{\underset{2\leq s \leq \sum_{i=1}^{r-2}p_{r,i}(\mu)+\ell_{s}}{\underset{\mu=\ell_1+\cdots+\ell_{s-\ell_s}}{m\leq \ell_{s} \leq \cdots \leq \ell_1 \leq n}}} \frac{q^{n^2+\ell_1+\cdots +\ell_s}}{(q)_n};$

  \begin{equation*}
\begin{aligned}
\sum_{\underset{1\leq j \leq \sum_{i=1}^{r-2}p_{r,i}(\mu)+l-1}{\underset{\mu=\ell_1+\cdots+\ell_{j-l+1}}{l\leq \ell_j \leq \cdots \leq \ell_1 \leq n}}} \frac{q^{n^2+\ell_1+\cdots +\ell_j+l}}{(q)_n}
&=\sum_{\underset{1\leq j \leq \sum_{i=1}^{r-2}p_{r,i}(\mu)+\ell_{j+1}-1}{\underset{\mu=\ell_1+\cdots+\ell_{j-\ell_{j+1}+1}}{m\leq \ell_{j+1} \leq \ell_j \leq \cdots \leq \ell_1 \leq n}}} \frac{q^{n^2+\ell_1+\cdots +\ell_j+\ell_{j+1}}}{(q)_n}\\
&=\sum_{\underset{2\leq s \leq \sum_{i=1}^{r-2}p_{r,i}(\mu)+\ell_{s}}{\underset{\mu=\ell_1+\cdots+\ell_{s-\ell_s}}{m\leq \ell_{s} \leq \cdots \leq \ell_1 \leq n}}} \frac{q^{n^2+\ell_1+\cdots +\ell_s}}{(q)_n}.
\end{aligned}
\end{equation*}

For each $s>\ell_s$ we have $\mu=\ell_1+\cdots+\ell_{s-\ell_s}.$ Let us denote the partition $\ell_1+\cdots+\ell_s$ by $\Lambda.$ By the definition of the $(r,i)-$new part of a partition we have:

$$p_{r+1,l}(\Lambda)=p_{r,l-1}(\mu),$$
where $1\leq l \leq r-1.$
%$p_{r+1,1}(\Lambda)=\ell_s,$

%$p_{r+1,2}(\Lambda)=\ell_{s-\ell_s}=p_{r,1}(\mu),$

%$p_{r+1,3}(\Lambda)=\ell_{s-\ell_s-\ell_{s-\ell_s}}=p_{r,2}(\mu),$

%$\vdots$

%$p_{r+1,r-1}(\Lambda)=p_{r,r-2}(\mu).$

So we can replace the last term of Equation \ref{in} by:

$$\sum_{\underset{2\leq s \leq \sum_{i=1}^{r-1}p_{r+1,i}(\Lambda)}{\underset{\Lambda=\ell_1+\cdots+\ell_s}{m\leq \ell_{s} \leq \cdots \leq \ell_1 \leq n}}} \frac{q^{n^2+\ell_1+\cdots +\ell_s}}{(q)_n}.$$
%So we have $\mathcal{C}_4=\sum_{\underset{2\leq s \leq \sum_{i=1}^{r-1}p_{r+1,i}(\Lambda)}{\underset{\Lambda=\ell_1+\cdots+\ell_s}{m\leq \ell_{s} \leq \cdots \leq \ell_1 \leq n}}} \frac{q^{n^2+\ell_1+\cdots +\ell_s}}{(q)_n},$ and $H_{r+1,1}^m=1+\sum_{i=1}^{4}\mathcal{C}_i,$ which is equal to:
So $H_{r+1,1}^m$ is equal to:

$$\Big(\sum_{n=0}^{m-1}\frac{q^{nm}}{(q)_n}+\sum_{n\geq m } \frac{q^{n^2}}{(q)_n}\Big)+(\sum_{m\leq \ell \leq n}\frac{q^{n^2+\ell}}{(q)_n})+\Big(\sum_{\underset{2\leq s \leq \sum_{i=1}^{r-1}p_{r+1,i}(\Lambda)}{\underset{\Lambda=\ell_1+\cdots+\ell_s}{m\leq \ell_{s} \leq \cdots \leq \ell_1 \leq n}}} \frac{q^{n^2+\ell_1+\cdots +\ell_s}}{(q)_n}\Big).$$

This proves the basic case of the induction on $r$. Assume now that the formula is true for $H_{r+1,c-1}^m$  and we prove it for $H_{r+1,c}^m.$ We have:

%$$H_{r+1,c}^m=1+q^mH_{r+1,c-1}+ q^{(m+1)}H_{r+1,c-1}^{m+1}+q^{m+2}H_{r+1,c-1}^{m+2}+\cdots$$

$$H_{r+1,c}^m=1+\sum_{l\geq m}H^l_{r+1,c-1}.$$
By the induction hypothesis on $c$ it is equal to:

\begin{equation}\label{H_{r+1}}
H_{r+1,c}^m=
1+\sum_{l\geq m}\sum_{n=0}^{l-1}\frac{q^{(n+1)l}}{(q)_n}+\sum_{l\geq m}\sum_{n\geq l}\frac{q^{n^2+l}}{(q)_n}+\sum_{l\geq m}\sum_{\underset{1\leq j \leq \sum_{i=1}^{r-1}p_{r+1,i}(\mu)+c-2}{\underset{\mu=\ell_1+\cdots+\ell_{j-c+2}}{l\leq \ell_j \leq \cdots \leq \ell_1 \leq n}}} \frac{q^{n^2+\ell_1+\cdots +\ell_j+l}}{(q)_n}.
\end{equation}

%$$\begin{matrix}

  %    1 +q^m\Big(1   & +\sum_{n=1}^{m-1}\frac{q^{nm}}{(q)_n} & +\sum_{n\geq m}\frac{q^{n^2}}{(q)_n}& +\sum_{\underset{1\leq j \leq \sum_{i=1}^{r-1}p_{r+1,i}(\mu)+c-2}{\underset{\mu=\ell_1+\cdots+\ell_{j-c+2}}{m\leq \ell_j \leq \cdots \leq \ell_1 \leq n}}} \frac{q^{n^2+\ell_1+\cdots +\ell_j}}{(q)_n} \Big)& \\
 %       +q^{m+1}\Big(1   & +\sum_{n=1}^{m}\frac{q^{n(m+1)}}{(q)_n} & +\sum_{n\geq m+1}\frac{q^{n^2}}{(q)_n}& +\sum_{\underset{1\leq j \leq \sum_{i=1}^{r-1}p_{r+1,i}(\mu)+c-2}{\underset{\mu=\ell_1+\cdots+\ell_{j-c+2}}{m+1\leq \ell_j \leq \cdots \leq \ell_1 \leq n}}} \frac{q^{n^2+\ell_1+\cdots +\ell_j}}{(q)_n} \Big)& \\
 %       +q^{m+2}\Big(1   & +\sum_{n=1}^{m+1}\frac{q^{n(m+2)}}{(q)_n} & +\sum_{n\geq m+2}\frac{q^{n^2}}{(q)_n}& +\sum_{\underset{1\leq j \leq \sum_{i=1}^{r-1}p_{r+1,i}(\mu)+c-2}{\underset{\mu=\ell_1+\cdots+\ell_{j-c+2}}{m+2\leq \ell_j \leq \cdots \leq \ell_1 \leq n}}} \frac{q^{n^2+\ell_1+\cdots +\ell_j}}{(q)_n} \Big)&\\            
 %      +~~~\cdots\\ 
%    \end{matrix}$$
%Once again we denote the sum of the $i$-th column by $\mathcal{C}_i.$ Note that $\mathcal{C}_1,\mathcal{C}_2$ and $\mathcal{C}_3$ are the same as $\mathcal{C}_1,\mathcal{C}_2$ and $\mathcal{C}_3$  in the proof of the basic case of the second induction and by the similar computation as before we obtain:
Note that the first three terms of Equation \ref{H_{r+1}} are the same as the first three terms of Equation \ref{in}. By a similar computations for its last term we obtain:

$$\sum_{l\geq m}\sum_{\underset{1\leq j \leq \sum_{i=1}^{r-1}p_{r+1,i}(\mu)+c-2}{\underset{\mu=\ell_1+\cdots+\ell_{j-c+2}}{l\leq \ell_j \leq \cdots \leq \ell_1 \leq n}}} \frac{q^{n^2+\ell_1+\cdots +\ell_j+l}}{(q)_n}
 =\sum_{\underset{2\leq s \leq \sum_{i=1}^{r-1}p_{r+1,i}(\mu)+c-1}{\underset{\mu=\ell_1+\cdots+\ell_{s-c+1}}{m\leq \ell_{s} \leq \cdots \leq \ell_1 \leq n}}} \frac{q^{n^2+\ell_1+\cdots +\ell_j+\ell_{s}}}{(q)_n}. $$
 
So we have $H_{r+1,c}^m$ is equal to:

%$$1+(\frac{q^m}{(q)_1})+\Big(\sum_{n=2}^{m-1}\frac{q^{nm}}{(q)_n}+\sum_{n\geq m } \frac{q^{n^2}}{(q)_n}\Big)+(\sum_{m\leq \ell \leq n}\frac{q^{n^2+\ell}}{(q)_n})+\Big( \sum_{\underset{2\leq s \leq \sum_{i=1}^{r-1}p_{r+1,i}(\mu)+c-1}{\underset{\mu=\ell_1+\cdots+\ell_{s-c+1}}{m\leq \ell_{s} \leq \cdots \leq \ell_1 \leq n}}} \frac{q^{n^2+\ell_1+\cdots +\ell_j+\ell_{s}}}{(q)_n}\Big)=$$

%$$1+\sum_{n=1}^{m-1}\frac{q^{nm}}{(q)_n}+\sum_{n\geq m}\frac{q^{n^2}}{(q)_n}+\sum_{\underset{1\leq j \leq \sum_{i=1}^{r-1}p_{r+1,i}(\mu)+c-1}{\underset{\mu=\ell_1+\cdots+\ell_{j-c+1}}{m\leq \ell_j \leq \cdots \leq \ell_1 \leq n}}} \frac{q^{n^2+\ell_1+\cdots +\ell_j}}{(q)_n}.$$
 
 \begin{equation*}
\begin{aligned}
\Big(\sum_{n=0}^{m-1}\frac{q^{nm}}{(q)_n}+\sum_{n\geq m } \frac{q^{n^2}}{(q)_n}\Big)+(\sum_{m\leq \ell \leq n}\frac{q^{n^2+\ell}}{(q)_n})\\
+\Big( \sum_{\underset{2\leq s \leq \sum_{i=1}^{r-1}p_{r+1,i}(\mu)+c-1}{\underset{\mu=\ell_1+\cdots+\ell_{s-c+1}}{m\leq \ell_{s} \leq \cdots \leq \ell_1 \leq n}}} \frac{q^{n^2+\ell_1+\cdots +\ell_j+\ell_{s}}}{(q)_n}\Big)\\
=\sum_{n=0}^{m-1}\frac{q^{nm}}{(q)_n}+\sum_{n\geq m}\frac{q^{n^2}}{(q)_n}+\sum_{\underset{1\leq j \leq \sum_{i=1}^{r-1}p_{r+1,i}(\mu)+c-1}{\underset{\mu=\ell_1+\cdots+\ell_{j-c+1}}{m\leq \ell_j \leq \cdots \leq \ell_1 \leq n}}} \frac{q^{n^2+\ell_1+\cdots +\ell_j}}{(q)_n}.
\end{aligned}
\end{equation*}    

\end{proof}

\begin{rem}\label{lead=}

On the one hand, in the second section we mentioned that:
 $$HP_{\frac{S}{I_{r,r}}}(q)=\sum_{n\geq 0}B_{r,r}(n)q^n.$$

On the other hand, in Example \ref{arc} we saw that $I_{r,r}$ is the leading ideal of the ideal $I_r$ with respect to the weighted reverse lexicographic order. So as we mentioned in the sixth section, we have $HP_{\frac{S}{I_r}}(q)=HP_{\frac{S}{I_{r,r}}}(q),$ and we guess that $I'_{r,i}$ is the leading ideal of $I_r$ with respect to the weighted lexicographic order. If we confirm this guess we will have:
%Using \textit{Singular}, we tried to find the form of the monomials which generate the leading ideal of
%$I_r,$ but this time with respect to the weighted lexicographic order. We guess that $I'_{r,r}$ is this new leading ideal. 
$$\sum_{n\geq0}C_{r,r}(n)q^n \underset{\text{By Proposition \ref{I'}}}{=}HP_{\frac{S}{I'_{r,r}}}(q)=HP_{\frac{S}{I_r}}(q)=HP_{\frac{S}{I_{r,r}}}(q)=\sum_{n\geq 0}B_{r,r}(n)q^n.$$ 
This would prove Conjecture \ref{conjG} for $i=r:$ 
\end{rem}

\begin{conj}\label{anal i=r} For $r\geq 3$ we have:
%$$H^1_{r,1}=1+\sum_{n\geq 1}\frac{q^{n^2}}{(q)_n}+\sum_{\underset{1\leq j \leq \sum_{i=1}^{r-2}p_{r,i}(\mu)}{\underset{\mu=\ell_1+\cdots+\ell_{j}}{1\leq \ell_j \leq \cdots \leq \ell_1 \leq n}}} \frac{q^{n^2+\ell_1+\cdots +\ell_j}}{(q)_n}=$$

%$$\sum_{n_1,n_2,\dots n_{r-1} \geq 0} \frac{q^{N_1^2+N_2^2+\dots +N_{r-1}^2}}{(q)_{n_1}(q)_{n_2}\dots (q)_{n_{r-1}}}=\prod_{\underset{n \not\equiv 0, \pm r (mod. 2r+1)}{n\geq1}} \frac{1}{1-q^n},$$

\begin{equation*}
\begin{aligned}
H^1_{r,1}& = 1+\sum_{n\geq 1}\frac{q^{n^2}}{(q)_n}+\sum_{\underset{1\leq j \leq \sum_{i=1}^{r-2}p_{r,i}(\mu)}{\underset{\mu=\ell_1+\cdots+\ell_{j}}{1\leq \ell_j \leq \cdots \leq \ell_1 \leq n}}} \frac{q^{n^2+\ell_1+\cdots +\ell_j}}{(q)_n} \\
      &=\sum_{n_1,n_2,\dots n_{r-1} \geq 0} \frac{q^{N_1^2+N_2^2+\dots +N_{r-1}^2}}{(q)_{n_1}(q)_{n_2}\dots (q)_{n_{r-1}}},\\
\end{aligned}
\end{equation*}

where $q$ is a variable, $N_j=n_j +n_{j+1} +\dots + n_{r-1}$ for all $1 \leq j\leq r-1 $ and $(q)_n=(1-q)(1-q^2)\cdots(1-q^n)$.

\end{conj}

%For the integers $4\leq r$ and $1\leq i \leq r,$ let $I_{r,i}$ be the following ideal:
%$$(x_1^i,x_1^{i-1}x_2^{r-i+1},\cdots,x_1x_2^{r-1},x_j^{r-\ell} x_{j+1}^{\ell}| \ 2\leq j, \ 0\leq \ell \leq r-1).$$
%Note that $I_{r,i}=(x_1^i,I_{r,r})$ and $I'_{r,i}=(x_1^i,I'_{r,r}).$

%Now let $I$ denote the ideal that defines the focussed arc algebra of $X$ defined by $(x^r)$ in $k[x].$
%In [BMS] (see also [BMS1]),  C. Bruschek, H. Mourtada, J. Schepers, proved that $I_{r,r}$ is the leading ideal of the ideal $I$ with respect to the weighted reverse lexicographic order. As we mentioned before in the second section, computing the Hilber-Poincar\'e series of the graded algebra $\frac{k[x_1,x_2,\cdots]}{I_{r,i}}$ is equivalent to counting the number of the partitions of  an integer $n\geq 0$ which are counted by $B_{r,i}(n)$.  

%We considered the weighted lexicographic order instead of the weighted reverse lexicographic order. Using \textit{Singular}, we tried to find the form of the monomials which generate the leading ideal of
%$I,$ but this time with respect to the weighted lexicographic order. We guess that $I'_{r,r}$ is this new leading ideal. If we confirm this guess we will have
%$$HP(\frac{k[x_1,x_2,\cdots]}{I_{r,r}})=HP(\frac{k[x_1,x_2,\cdots]}{I})=HP(\frac{k[x_1,x_2,\cdots]}{I'_{r,r}}).$$ 
%Which proves Conjecture \ref{conjG} for $i=r.$

\hfill \break

\noindent Equipe G\'eom\'etrie et Dynamique, \\
Institut Math\'ematique de Jussieu-Paris Rive Gauche,\\
 Universit\'e Paris Diderot, \\
 B\^atiment Sophie Germain, case 7012,\\
75205 Paris Cedex 13, France.\\

\noindent Email : pooneh.afshari@gmail.com
 
\end{document}